\documentclass[12pt,a4paper]{article}

\usepackage{a4,amsmath,amsfonts,amssymb,epsfig}
\usepackage[latin1]{inputenc}

\newcommand{\RR}{{\mathbb{R}}}
\newcommand{\CC}{{\mathbb{C}}}

\newcommand{\Ss}{{\mathbb{S}}}

\newcommand{\proof}{{\noindent\bf Proof: }}

\def\qed{\unskip\nobreak\hfil\penalty50\hskip1.75em\null\nobreak\hfil
$\blacksquare$ {\parfillskip=0pt \finalhyphendemerits=0 \par}\medbreak}

\newcommand\capsize{\relax}
\newcommand\nobf{\noindent\bf}

\newcommand\sgn{\operatorname{sgn}}
\newcommand\grad{\operatorname{grad}}
\newcommand\interior{\operatorname{int}}

\newtheorem{theo}{Theorem}
\newtheorem{lemma}[theo]{Lemma}
\newtheorem{prop}[theo]{Proposition}

\title{Functions from $\RR^2$ to $\RR^2$: \goodbreak a study in nonlinearity}
\author{Nicolau C. Saldanha and Carlos Tomei}

\begin{document}
\maketitle

% \begin{abstract}
% Rouch\'e's theorem and winding numbers are standard tools for
% counting solutions of holomorphic equations in the complex plane.
% We explore counterparts to these tools for smooth systems of two real
% equations in two real unknowns.
% In the process, we outline the implementation of a
% system solver which makes use of this topological information.
% \end{abstract}
%
% MSC 26B99, 58K15, 65H10, 65H20, 90C30.
% Keywords:
% Rouch\'e's theorem, cusp, fold, zeros of functions, continuation methods

\section{Introduction}

Calculus students learn how to draw graphs of functions from
$\RR$ to $\RR$ and undergraduates studying complex variable
learn about geometric properties of functions like $f(z) = z^3$
and $g(z) = e^z$. Some teachers go further and introduce
a few examples of conformal mappings. A picture is worth a
thousand words, but more can be said on their favor: they provide
a good exercise in combining theoretical facts in a consistent
fashion. Indeed, to obtain the graph of a real function, a student
considers its derivatives, asymptotic behavior and some
special points, among other features. Something similar happens
in the study of conformal mappings.

In this text, we consider functions from $\RR^2$ to $\RR^2$ and
along the way assemble a number of tools from undergraduate courses.
We describe a graphical representation of such functions and,
for functions which are visually too complicated,
we still count preimages, in a manner reminiscent of Rouch\'e 's theorem.
Why is it that such aspects of functions from the plane to the plane
are not more familiar? A reason might be the following. Most of the
information we compute about functions from the line to itself,
or about holomorphic functions, concerns special points---typically
critical points, where the derivative is zero.
In the case of functions from the plane to the plane,
we need to consider critical curves,
where the Jacobian matrix is not invertible.
Such curves are often impossible to describe in simple closed form.

Enter the computer: we should think of the study of a given function
from the plane to the plane as a description of certain relevant
objects, in a way that these objects become amenable to numerics.
In this sense, the time is ripe for this new case
study in nonlinear theory, in the same way that we feel more at
ease nowadays with showing students how to evaluate roots of
polynomials of degree $6$, or eigenvalues of $5 \times 5$ matrices.

The theory should operate on two levels: we should learn
enough to get qualitative information about simple examples, and
we should be able to derive numerical procedures to handle general
cases. In particular, such procedures should extend our knowledge
of the preimages of a point, from mere counting to explicit
computation.

In section 2 we present a representative function $F_0$
which will be our favorite test case throughout the paper;
in section 10 some additional examples are discussed.
Some of the tools required for this project belong to
the standard curriculum, others are just ahead.
All of them are basic when dealing with nonlinear problems.
Thus, for example, in section 3,
we describe the local behavior of a function at {\it folds} and {\it cusps},
special critical points in the domain where the inverse
function theorem does not apply.  We will compute winding numbers and 
will also consider, in section 8, the {\it rotation number} of a $C^1$ curve.
Some aspects of covering space theory,
presented in sections 5 and 6,
will help us fit together local information.
In particular, we will be able to perform {\it compatibility checks},
discussed in section 9, which often indicate the presence of
yet unknown critical curves.

Some theoretical aspects have computational counterparts.
For example, under appropriate hypothesis,
the inverse function theorem asserts that a function is locally invertible
while Newton's method may be used to actually perform the inversion.
More generally, the implicit function theorem verifies
the regularity of critical curves and a predictor-corrector method
then traces the curve, as in section 4.
Similarly, covering space theory is closely related to
numerical continuation methods, employed in section 7.
Due to space limitations, we handle numerical aspects rather superficially,
providing sketches of arguments and indicating more specific literature.
Some results are quoted from standard references but we present
proofs of a few statements which are harder to find in book form.
Senior undergraduates should be able to follow through the
arguments.

Together with Iaci Malta, the authors have published
more technical texts (\cite{MST1}, \cite{MST2}).
The program (in rough form) which generated
pictures and computations for this paper is available (\cite{2x2}).
Both theoretical and computational aspects can be extended
to the study of functions from a bounded subset of the plane
to the plane (\cite{Duczmal}).
For a more general study of the geometry of functions between 
two surfaces, see \cite{FT}.

\section{A first example}

\begin{figure}[ht]
\vglue 11pt
\begin{center}
\epsfig{height=32mm,file=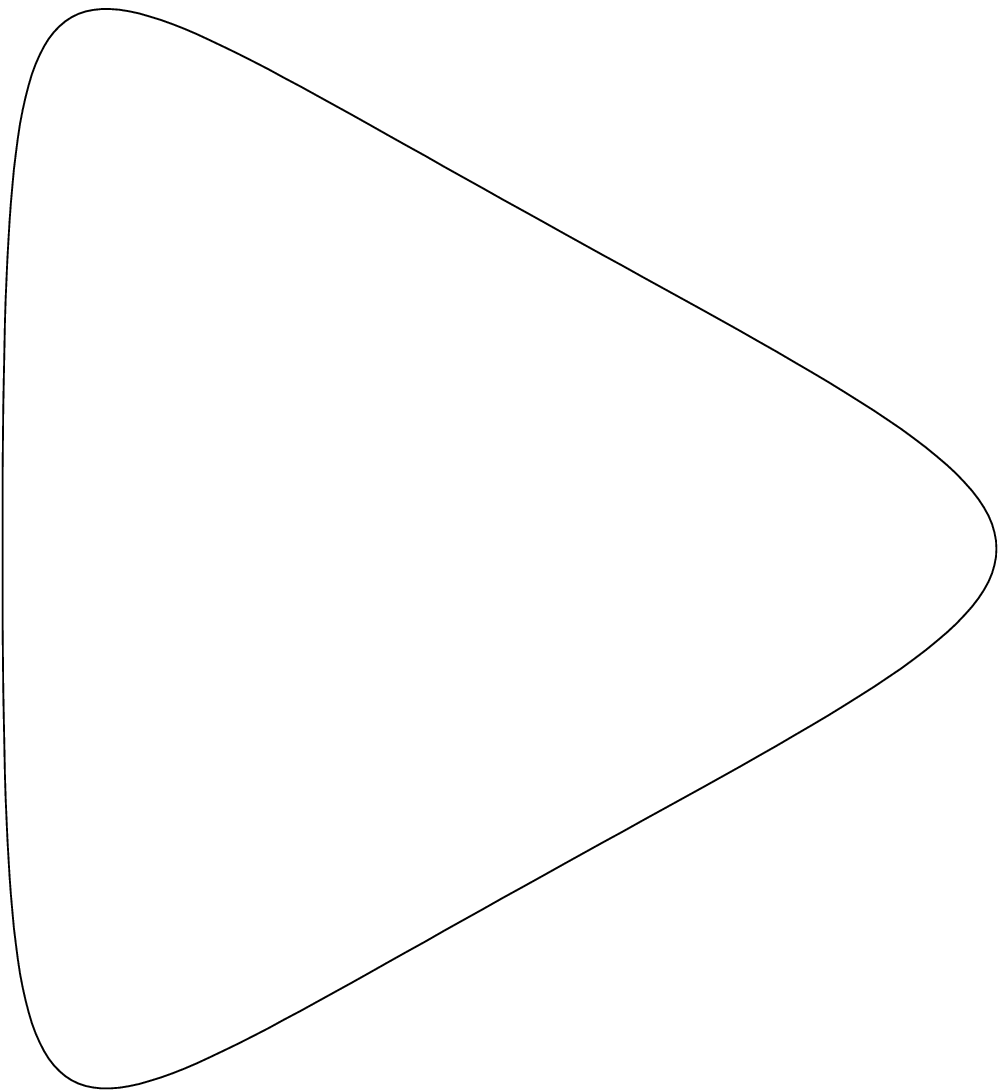} \qquad
\epsfig{height=32mm,file=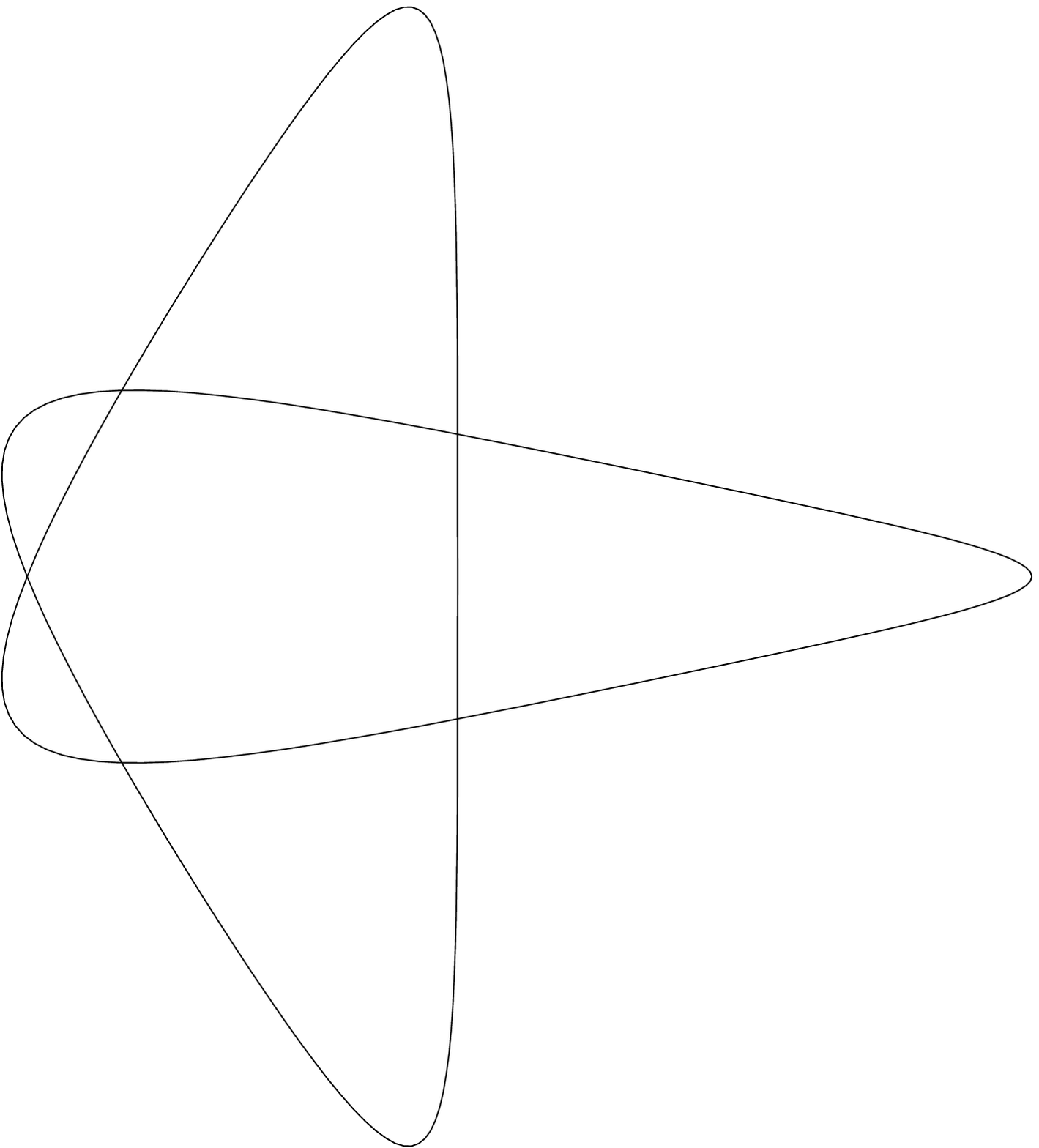} \qquad
\epsfig{height=32mm,file=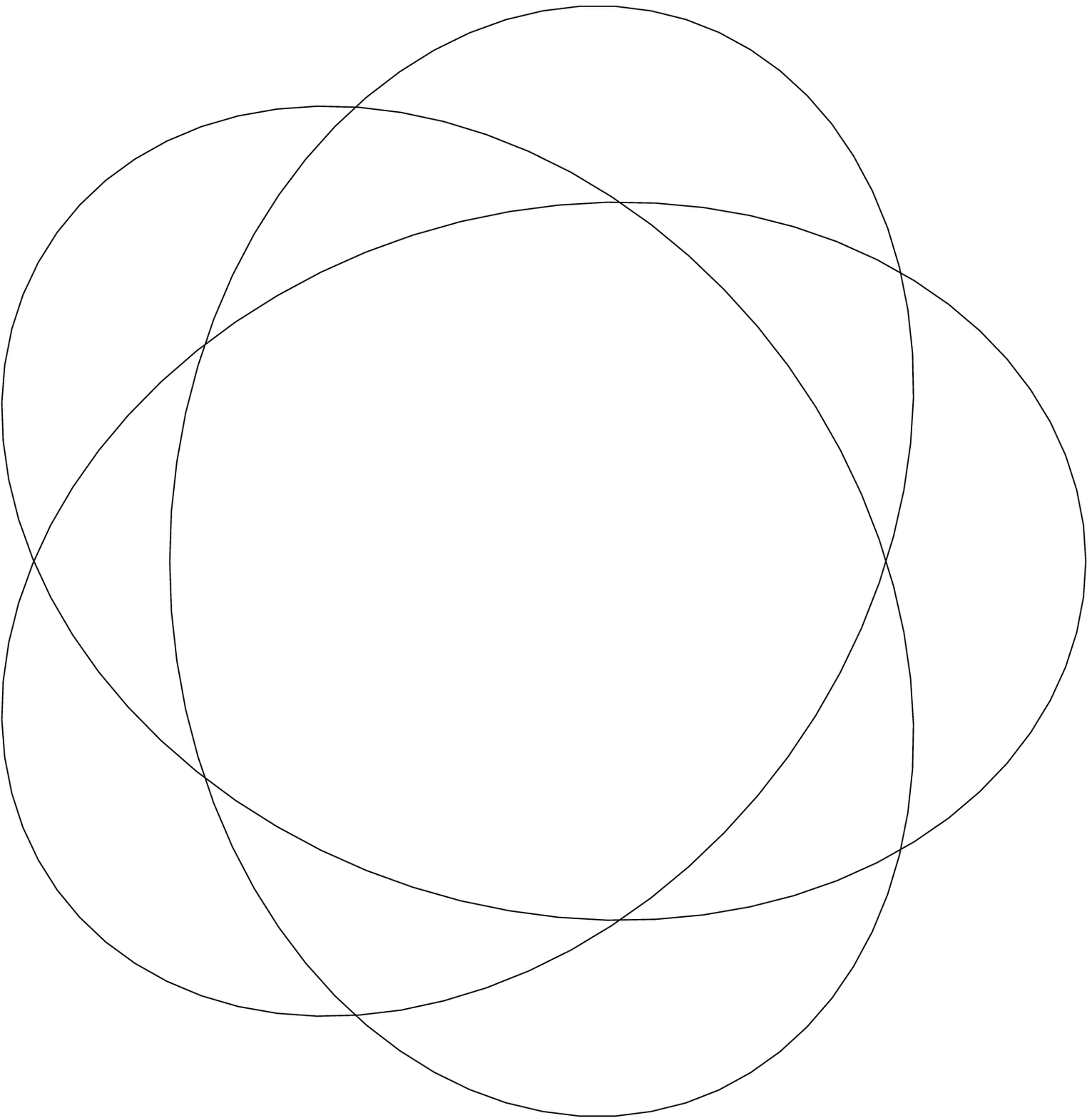}
\end{center}
\caption{\capsize Three different regimes (radii $0.1$, $1$ and $10$)}
\label{fig:3regs}
\end{figure}

Our first and favorite example is the function
\[ \begin{matrix}
F_0: & \RR^2 & \to & \RR^2 \\
& \begin{pmatrix}x \\ y \end{pmatrix} & \mapsto &
\begin{pmatrix}x^3 - 3 x y^2 + 2.5 x^2 - 2.5 y^2 + x \\
3 x^2 y - y^3 - 5 xy + y \end{pmatrix} \end{matrix} \] which,
in complex notation, can be written as
$F_0(z) = z^3 + 2.5 \bar{z}^2 + z$.
Due to the presence of $\bar{z}$, $F_0$ is not holomorphic.
As every Rouch\'e fan would notice,
$F_0$ acts on concentric circles centered at the origin
according to (at least) three different regimes: figure
\ref{fig:3regs} shows the images of circles with radii equal to
$0.1$, $1$ and $10$, respectively (the figures are not in scale).

% fcirc015 fcirc1 fcirc3

Indeed, orient the circles in the domain positively (i.e.,
counterclockwise). For radii close to $0$, $F_0$ takes these
circles to simple closed curves with winding number $1$ with
respect to $0$: this is clear from the fact that, close to the
origin, $F_0$ is essentially the identity. For radii close to $1$,
the term $2.5 \bar{z}^2$ dominates the other two ($|2.5 \bar{z}^2|
> |z^3| + |z|$ for $|z| \approx 1$) and one should expect the
images of such circles to be closed curves winding twice around
the origin with negative orientation---the winding number with
respect to $0$ of these curves is $-2$. Finally, $F_0$ takes
circles of large radius to curves winding three times positively
around the origin: $F_0$ near infinity looks like $z^3$.

% fcirc03 04 06 13 15 18
\begin{figure}[ht]
\vglue 11pt
\begin{center}
\epsfig{height=32mm,file=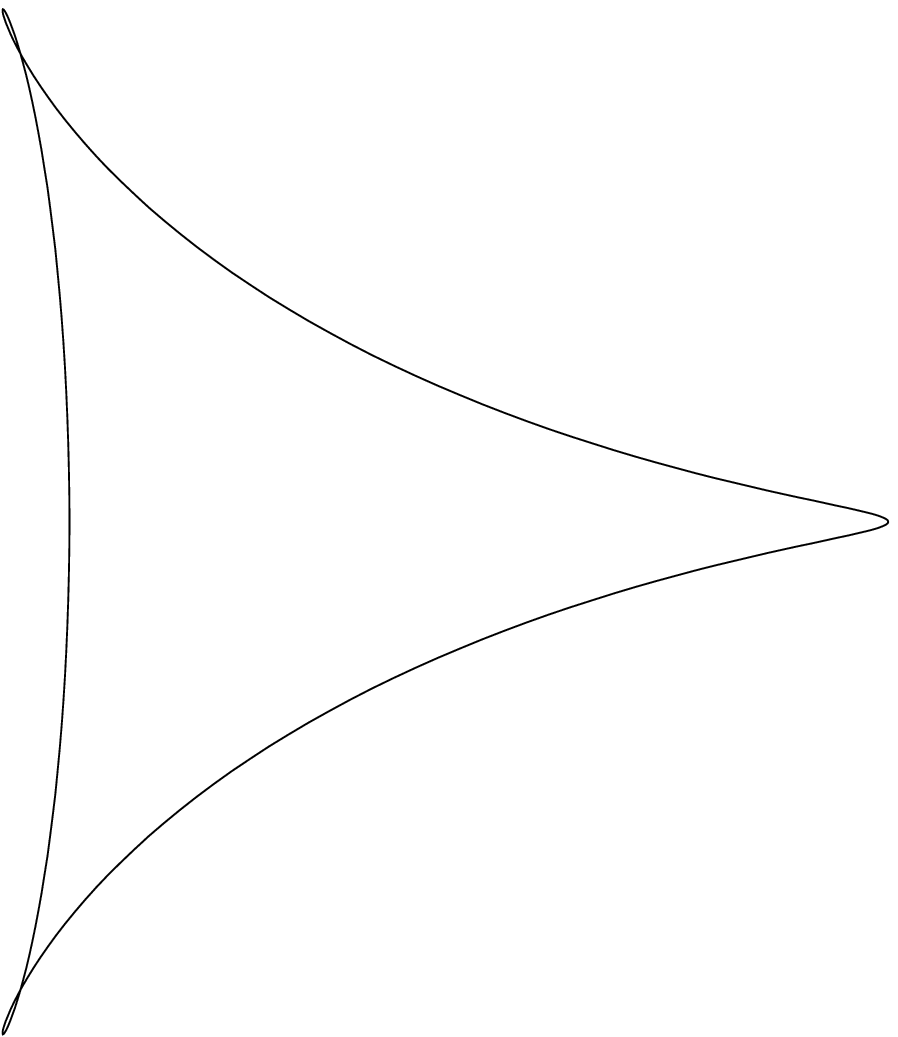} \qquad
\epsfig{height=32mm,file=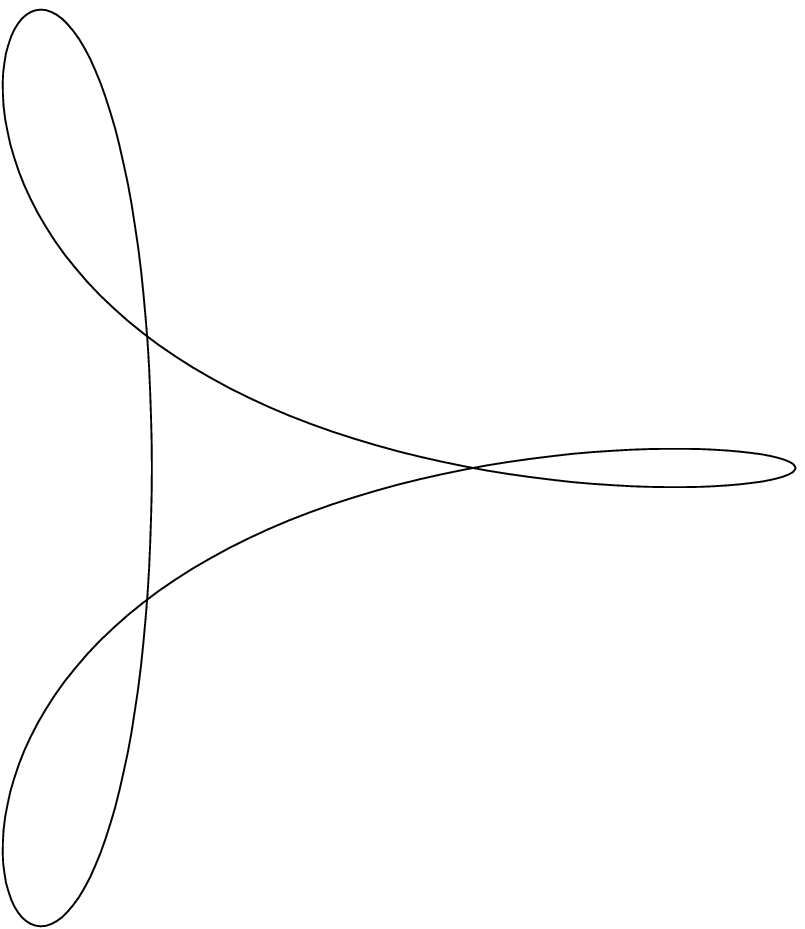} \qquad
\epsfig{height=32mm,file=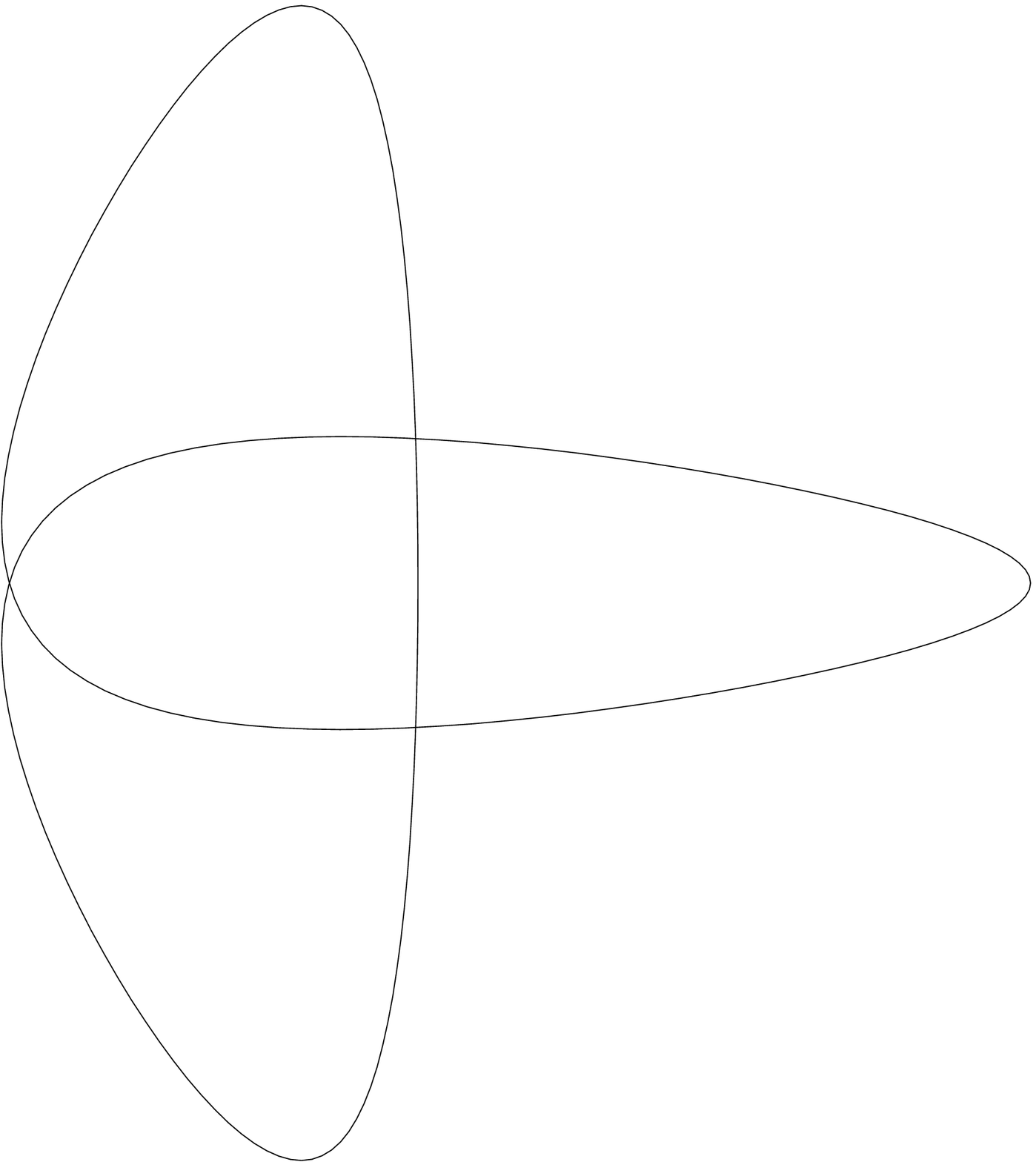}
\end{center}
\begin{center}
\epsfig{height=32mm,file=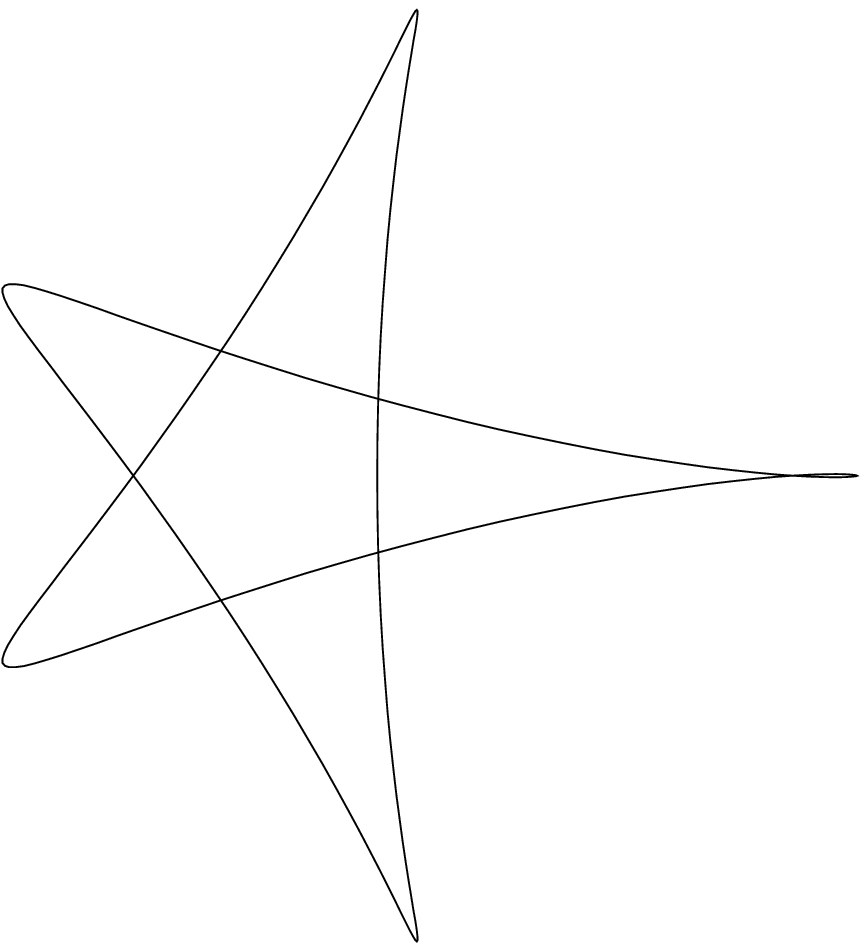} \qquad
\epsfig{height=32mm,file=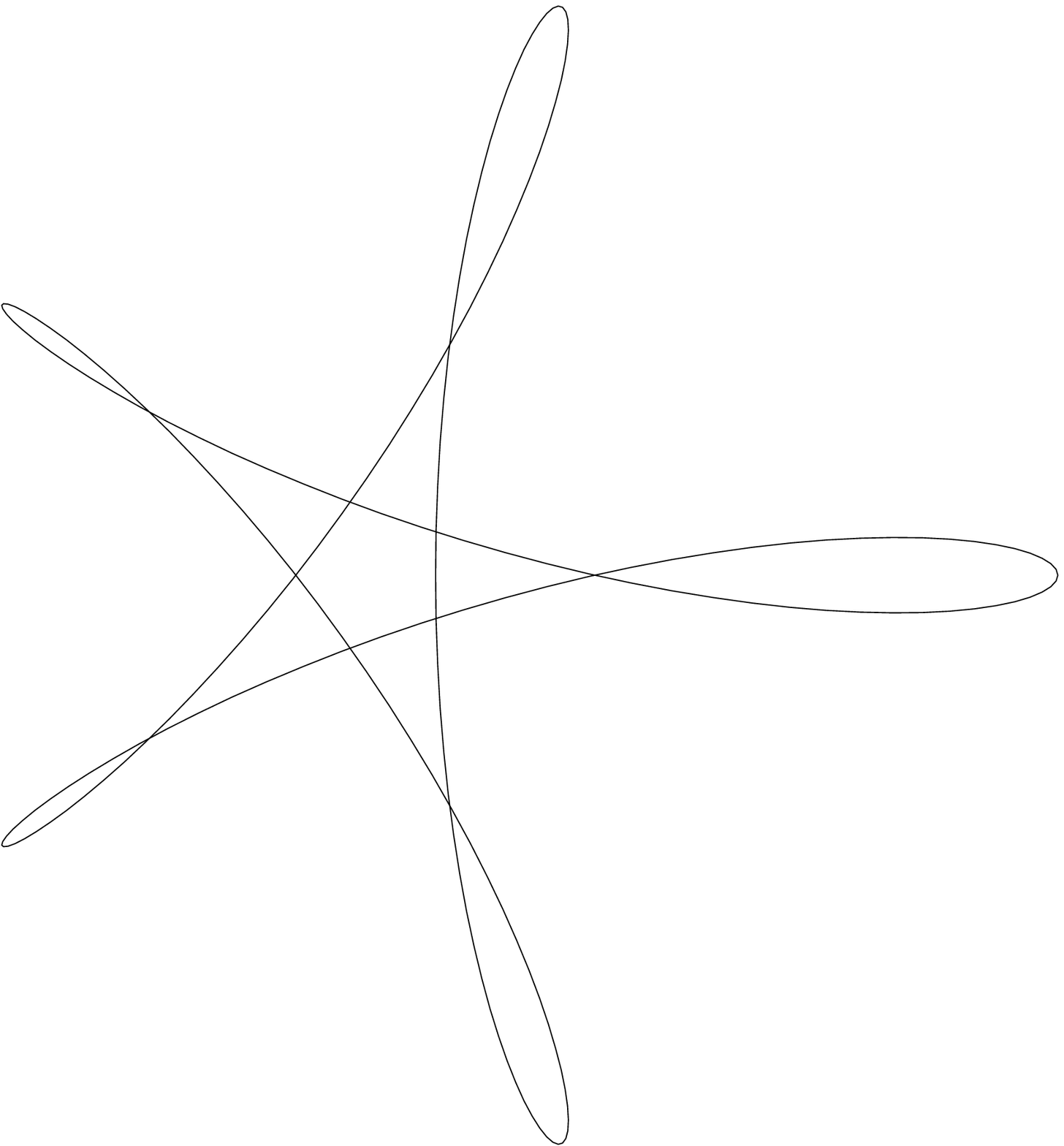}
\end{center}
\caption{\capsize Interpolating regimes:
radii $0.2$, $0.3$, $0.7$, $1.5$ and $2$ }
\label{fig:5regs}
\end{figure}

How are the images of concentric circles changing from one regime
to another? Figure \ref{fig:5regs} shows the images of
intermediate circles of radii $0.2$, $0.3$, $0.7$, $1.5$ and $2$.
It may be hard to see in the picture, but there are two small loops
in the first figure and one in the fourth; the five curves are
indeed smooth. We will address the formation of these patterns in
the sequel. We will draw other, more informative pictures and
will infer, for example, that the equation $F_0(z) = 0$ has nine solutions.

\section{Local theory}

We will use some basic facts of the local theory of functions in
the plane, obtained by Whitney in 1955 (\cite{W}).
The subject developed considerably under the name of singularity
theory after the work of Thom and Mather in the sixties.
An excellent reference with emphasis in applications is \cite{GS};
a more technical one is \cite{GG}.

Recall that, for a smooth function $F:\RR^2 \to \RR^2$, a point
$p$ is {\it regular} if the Jacobian $DF(p)$ is an invertible matrix.
From the inverse function theorem (\cite{L}, chap. XVII, \S 3, pg. 349),
after smooth changes of
variable in appropriate neighborhoods of a regular point $p$ and
its image $F(p)$, the function $F$ takes the form $\tilde{F}(x,y) = (x,y)$;
more precisely, there exist local diffeomorphisms $\Phi$ and $\Psi$
with $F = \Phi \circ \tilde{F} \circ \Psi$ as above.
Points which are not regular are {\it critical} and they
form the {\it critical set} $C$. A critical point $p_f$ is a {\it fold
point} (or, more informally, a fold) if, after changing variables
near $p_f$ and $F(p_f)$, $F$ becomes $\tilde{F}(x,y) = (x,y^2)$. Also,
a critical point $p_c$ is a {\it cusp point} (or, again, simply a
cusp) if changes of variables convert $F$ into $\tilde{F}(x,y) =
(x, y^3 - xy)$. The formulae for $\tilde{F}$ are the {\it normal forms}
of a function at a fold and at a cusp:
they imply that in appropriate neighborhoods of folds and cusps,
the critical set is a smooth arc consisting of folds and cusps;
also, cusps are isolated.

In the same
way that the hypothesis of the inverse function theorem guarantees
a simple normal form of a function near a regular point, there
are explicit conditions which characterize folds and cusps.
For example, a critical point $p_f$ is a fold of a
smooth function $F:\RR^2 \to \RR^2$ if two conditions hold. First,
the gradient of $\det DF$ should be nonzero at $p_f$, which implies,
from the implicit function theorem, that the critical set near $p_f$
is a curve. Second, $\grad\det DF(p_f)$ should not be orthogonal to
$\ker DF(p_f)$. There is a similar, more complicated
characterization of cusps, which we omit.

\begin{figure}[ht]
\vglue 11pt
\begin{center}
\epsfig{height=30mm,file=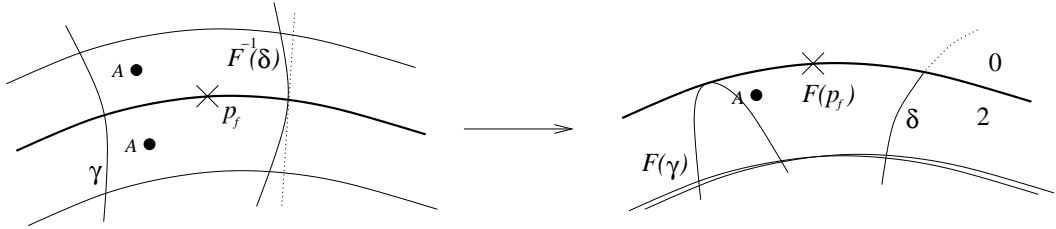}
\end{center}
\caption{\capsize Local behavior near a fold}
\label{fig:fold}
\end{figure}

A function near a fold behaves in a simple way.
All properties described below can be checked by referring to
the normal form. In
figure \ref{fig:fold}, a small arc of the critical set and
its image under the function $F$ are indicated with thick lines
for both types of critical points. Points in the domain with the
same image are indicated by the same label. The thinner lines on
both sides of the critical curves are taken to thin lines as
indicated. Near
a fold $p_f$, the function $F$ takes points to a single side of the
image of the critical arc. Thus, a point $w$ near $F(p_f)$ has $0$,
$1$ or $2$ preimages near $p_f$, depending on its position with
respect to the image of the critical arc.
The image of a curve $\gamma$ {\it transversal} to the critical set
at a fold point is generically a nonsingular curve $F(\gamma)$
tangent to $F(C)$
(two smooth arcs are transversal at an intersection point
if their tangent vectors are linearly independent).
The inverse image of a curve $\delta$ transversal to $F(C)$
is a curve tangent to $\ker(DF)$ at $C$.
Only one side of $\delta$ actually has preimages
(in the figure, the dotted part of $\delta$ is not in the
image of $F$ near $p_f$).

\begin{figure}[ht]
\vglue 11pt
\begin{center}
\epsfig{height=30mm,file=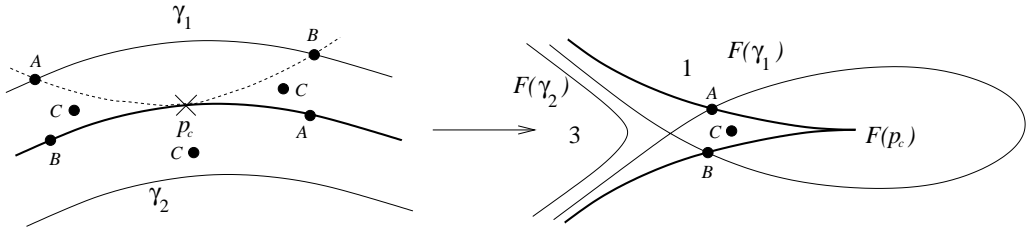}
\end{center}
\caption{\capsize Local behavior near a cusp}
\label{fig:cusp}
\end{figure}

Points $w$ near the
image $F(p_c)$ of a cusp $p_c$ may have $1$, $2$ or $3$ preimages near
$p_c$. Arcs $\gamma_1$ and $\gamma_2$ in figure \ref{fig:cusp}
have qualitatively different images:
$F(\gamma_1)$ undergoes a loop around $F(p_c)$, $F(\gamma_2)$
does not. We also indicate the (nontrivial)
preimage of the image of the critical curve near the cusp: notice
that it lies to one side of the critical curve. We say that the
cusp is {\it effective} on that side.

Whitney defined {\it excellent functions} as functions having only
folds and cusps as critical points. The regular points of an
excellent function $F: \RR^2 \to \RR^2$ form an open dense subset
of the plane and its critical set $C$ is a disjoint union of
isolated smooth curves. It is easy to see from the normal form
that cusps form a discrete set and the rest of $C$ consists of
arcs of folds. Fortunately, excellent functions from the plane to
the plane are abundant: this allows us to ignore more complicated
critical points.

\begin{theo}[Whitney,\cite{W}]
In the $C^r$ topology on compact sets ($r \ge 3$), the set
of excellent functions $F: \RR^2 \to \RR^2$ is residual.
\end{theo}

It is instructive to compare the local theory of excellent
functions to the local theory of holomorphic functions: we remind
the reader of the normal form of a holomorphic function at a
critical point.

\begin{prop}
Let $f: A \to \CC$ be a holomorphic function with a critical
point $z_0$ which is a zero of order $n-1$ of $f'$.
Then its normal form is given by $\tilde{f}(w) = w^n$.
More precisely, there exist local holomorphic diffeomorphisms
$\phi$, $\psi$ with $\phi(0) = f(z_0)$, $\psi(z_0) = 0$
for which $f = \phi \circ \tilde{f} \circ \psi$
in a neighborhood of $z_0$.
\end{prop}

The open mapping theorem and the maximum modulus theorem
(\cite{BN}, chapters 6 and 7) follow easily from this local form.

{\nobf Proof: }
Write
\[ f(z) = a + b (z-z_0)^n g(z) = a + b \big ((z-z_0) h(z) \big)^n, \quad
a = f(z_0), \; b = f^{(n)}(z_0)/n!, \]
where $g$ and $h$ are
holomorphic functions with $g(z_0) = h(z_0) = 1$ and $g(z) = (h(z))^n$.
Now set $\phi(u) = a+bu$, $\psi(z) = (z-z_0) h(z)$ and we are done.
\qed

In particular, critical points of nonconstant holomorphic
functions are isolated and they certainly may not be folds or
cusps. How can Whitney´s theorem be true then? An excellent
function $F$ near a holomorphic function $f$ must be
non-holomorphic.
For instance, to approximate $f(z) = z^7$ by an
excellent function, one may try $F(z) = z^7 + \epsilon \bar{z}$,
which indeed works for small $\epsilon$. There is a natural 
counterpart to Whitney's theorem for holomorphic functions:
in a residual set of holomorphic functions, the second derivative
is nonzero at all critical points.

\section{Tracing the critical set}

\begin{figure}[ht]
\vglue 11pt
\begin{center}
\epsfig{height=16mm,file=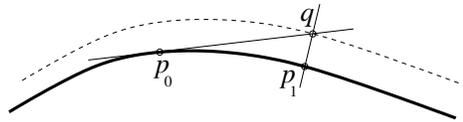}
\end{center}
\caption{\capsize A simple predictor-corrector method}
\label{fig:pc}
\end{figure}

Searching for critical curves by hand is hard even
for a polynomial map of low degree, such as our $F_0$.
Classifying critical points as folds,
cusps or yet something else is even harder.
A more practical approach is to go through numerical computations.
% to check,
% up to a degree of satisfaction, that critical points satisfy the
% characterizations of folds and cusps mentioned above.
In order to study the critical set of a function $F$,
our program first searches for points $p_+$ and $p_-$
for which $DF(p_+)$ and $DF(p_-)$ have determinants of opposite sign.
By continuity, there must be a critical point $p_0$ (i.e., $\det DF(p_0)=0$)
in the segment joining $p_+$ and $p_-$. After computing $p_0$,
the program obtains some points $p_1, p_2, \ldots$ in the
critical curve through $p_0$
(i.e., the level through $p_0$ of $\det DF_0$)
by a {\it predictor-corrector method}
(a fine presentation of this class of methods is given in \cite{AG}).
A very simple example of this technique is the following.
As in figure \ref{fig:pc}, draw a tangent line to the
critical curve through $p_0$ and take a point $q$ on this line at
a short distance $h$ from $p_0$. Now through $q$ draw a second
line parallel to $\grad \det DF(q)$ and on this line solve
$\det DF(p_1) = 0$ by Newton's method with initial condition $q$.

As the critical points $p_0, p_1, \ldots$ are computed, the
program checks the conditions characterizing folds
for segments from $p_i$ to $p_{i+1}$:
this is done in order to detect cusps (or other singularities).
On segments which do not pass this test,
the program searches for cusps and validates them with
additional tests which we do not detail.
These tests ascertain
with considerable robustness and reliability that all critical
points on this critical curve are indeed folds or cusps.

\begin{figure}[ht]
\vglue 11pt
\begin{center}
\epsfig{height=50mm,file=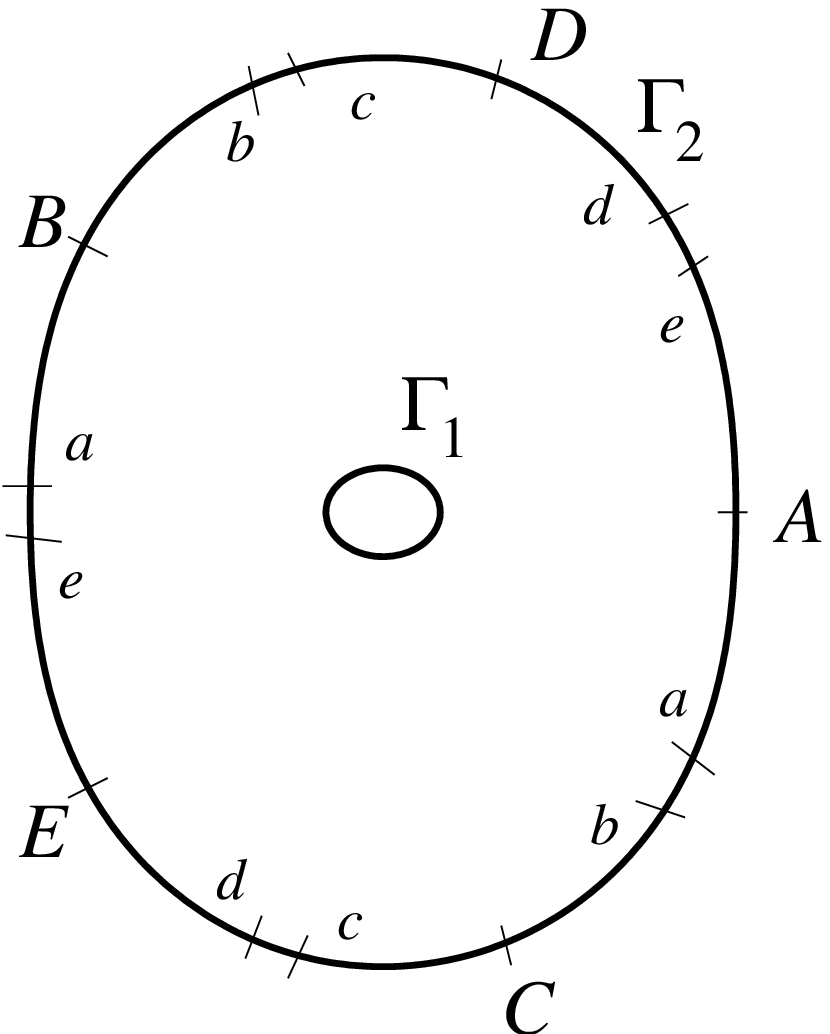}
\qquad \qquad
\epsfig{height=50mm,file=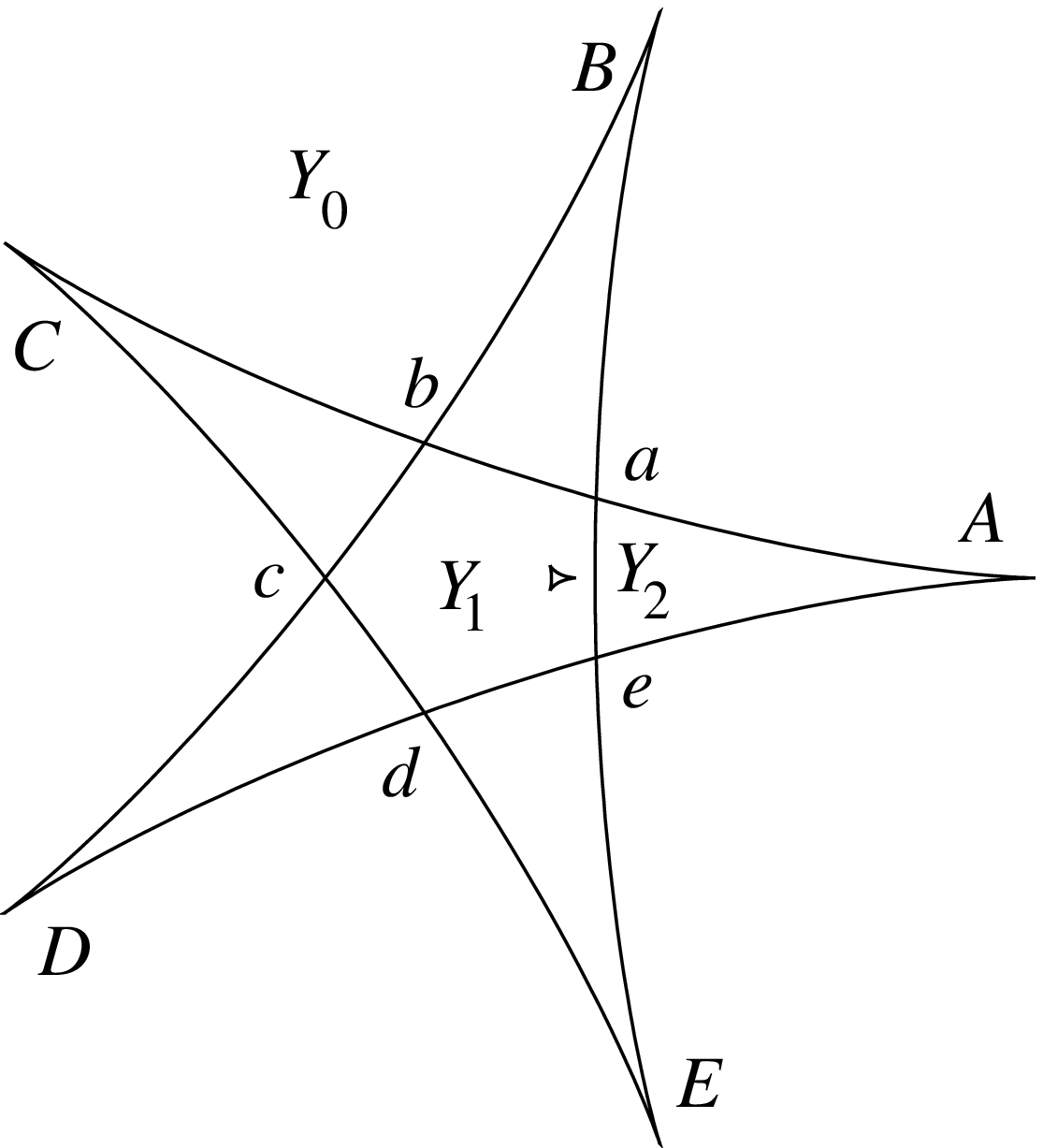}
\end{center}
\caption{\capsize The critical curves of $F_0$ and their images}
\label{fig:criticaldom}
\end{figure}

What is the critical set of $F_0$?
And for that matter, is it even excellent?
The left part of figure \ref{fig:criticaldom} shows both critical
curves $\Gamma_1$ and $\Gamma_2$ of the function $F_0$: they are
ovals around the origin. The images of the critical curves
are on the right: $F_0(\Gamma_1)$ is a small curvilinear triangle
surrounding the origin and $F_0(\Gamma_2)$ is a stellated
pentagon. Indeed, numerics confirm that the two critical curves
have 3 and 5 cusps. This is in agreement with the almost polygonal
shape of the images of circles of radii $0.2$ and $1.6$ in figure
\ref{fig:5regs}. The labels on the outer critical curve $\Gamma_2$
are of two kinds: capital letters indicate cusps and lower case
letters are preimages of self-intersections of $F_0(\Gamma_2)$.
The three cusps on $\Gamma_1$ are not indicated, but the reader
may check that if $\Gamma_1$ is traversed counterclockwise then so
is $F_0(\Gamma_1)$.

Critical curves are thus described by lists of points, some of them cusps.
Outside the known part of the critical set,
however, very few points have been considered.
Without a labor intensive search, how do we know
if we have found all critical curves of a function?
This, in general, is a nontrivial issue and we shall say
more about it in section 9.

\section{Counting preimages}

We set the
information obtained so far in a more robust setting
lest the reader think that we are letting pictures
take control over mathematical reasoning.
We begin by stating without proof a stronger form of
the Jordan curve theorem.

\begin{theo}
\label{theo:schoenflies}
Let $\gamma \subset \RR^2$ be a simple closed curve.
The curve $\gamma$ is the boundary of a closed topological disk $D$.
The open set $\RR^2 - \gamma$ has precisely two connected components:
the interior of $D$ and the complement of $D$.
Furthermore, if $\gamma$ is a piecewise smooth curve then
there exists a homeomorphism from $D$ to the closed unit disk
whose restriction to the interior of $D$ is holomorphic.
\end{theo}

The first claim is known as the Schoenflies theorem
for which a nice proof is given in \cite{Thomassen}.
The second is the standard Jordan theorem
(theorem 13.4, chapter 8, \cite{Munkres})
and the third is an extension of the Riemann mapping theorem
(14.19 in \cite{Rudin}).

We denote by $D_{\gamma}$ the closed (topological) disk surrounded
by the simple closed curve $\gamma$ and by $\interior D_{\gamma}$
the corresponding open disk.
The lemmas that follow are standard,
but somewhat hard to pinpoint in the literature.

\begin{lemma}
\label{lemma:homeo} Let $\gamma$ be a smooth simple closed curve
in $\RR^2$. Let $F: D_{\gamma} \to \RR^2$ be a $C^0$ map which is
$C^1$ in $\interior D_\gamma$. Assume that $F$ has no critical
points in $\interior D_\gamma$ and is injective on $\gamma$. Then
$F$ is a homeomorphism from $D_{\gamma}$ to its image $F(D)$.
\end{lemma}

\proof For readers acquainted with degree theory, the proof is
simpler; we sketch a more elementary argument.
Set $\delta = F(\gamma)$: clearly, $\delta$ is a simple, closed curve,
surrounding a closed (topological) disk $D_{\delta}$.
At every point $p \in \interior D_{\gamma}$, $F$ is open,
i.e., a small open ball around $p$ is
taken bijectively to a small open set around $F(p)$: this follows
from the inverse function theorem, since $F$ has no critical
points in $\interior D_{\gamma}$. Thus, any boundary point of
$F(D_{\gamma})$ ought to be in $F(\gamma) = \delta$.
By theorem \ref{theo:schoenflies}, the compact set $F(D_{\gamma})$
equals either $\delta$ or $D_\delta$; on the other hand,
$F(D_\gamma) = \delta$ is impossible, since the interior of $\delta$ is empty.

We now prove that $F$ is injective: from the
arguments above and the injectivity on $\gamma$, we only have to
show that if $p_0, p_1 \in \interior D_{\gamma}$ are such that $F(p_0) =
F(p_1)$ then $p_0 = p_1$. Let $\zeta: [0,1] \to \interior D_\gamma$ be a
smooth path with $\zeta(0) = p_0$, $\zeta(1) = p_1$ so that $(F \circ
\zeta)(0) = (F \circ \zeta)(1)$. Let $H: [0,1]^2 \to \interior
D_\Gamma$ be a smooth function with $H(0,t) = (F \circ \zeta)(t)$,
$H(s,0) = H(s,1) = H(1,t)$: the existence of such $H$ is ascertained
by theorem \ref{theo:schoenflies}.
We now construct $\zeta_s: [0,1] \to \interior D_\gamma$
so that $F(\zeta_s(t)) = H(s,t)$:
$\zeta_s$ is the solution of the differential equation
\[ \zeta_s(0) = p_0, \quad
\zeta'_s(t) =
(DF(\zeta_s(t)))^{-1} \frac{\partial H}{\partial t}(s,t). \]
Now $\zeta_1$ is
constant whence $\zeta_1(1) = p_0$. But $\zeta_s(1)$ depends
continuously on $s$ and satisfies $F(\zeta_s(1)) = F(p_0)$ for all $s$.
Therefore $\zeta_s(1) = p_0$ for all $s$ and $p_1 = \zeta_0(1) = p_0$.

Since $F$ is a continuous bijection from the compact set
$D_\gamma$ to the Hausdorff space $D_\delta$, $F$ is a
homeomorphism (theorem 5.6, chapter 3, \cite{Munkres}). \qed

A continuous function $F: \RR^2 \to \RR^2$ is \textit{proper} if
the inverse of any compact set is compact. The reader should have
no difficulty in proving that this is equivalent to saying that
$\lim_{p \to \infty} F(p) = \infty$. Our function $F_0$ is proper.

\begin{lemma}
\label{lemma:finite} Let $F: \RR^2 \to \RR^2$ be a proper
excellent function. The number of preimages under $F$ of any point
of $\RR^2$ is finite.
\end{lemma}

\proof By properness, all preimages of a point $w$ belong to a
closed disk $D$. If there are infinitely many of them, they must
accumulate at a point $p$, which, by continuity of $F$, is also a
preimage of $w$. Now, $p$ may not be either regular, a fold or a
cusp, since the three normal forms do not allow for infinitely
many local preimages, in disagreement with the excellence of $F$.
\qed

Given a closed set $X \subset \RR^2$,
we call the connected components of $\RR^2 - X$ the {\it tiles for} $X$.

\begin{lemma}
\label{lemma:count} Let $F: \RR^2 \to \RR^2$ be a proper excellent
function with critical set $C$.  On each tile $A$ for $F(C)$, the
number of preimages of $F$ is a constant.
\end{lemma}

\proof By connectivity, it suffices to show that the number of
preimages of points near $w \in A$ is constant. Let $p_1, \ldots,
p_k$ be the (finitely many) preimages of $w$. By hypothesis, they
are regular points, and thus there are open disjoint neighborhoods
$V_i, i=1,\ldots,k$, with $p_i \in V_i$ and so that $F$ restricts
as a homeomorphism from each $V_i$ to an open neighborhood $W$ of
$w$. Thus, points in $W$ have at least as many preimages as $w$.
Suppose now that for a sequence $w_j \in W$ converging to $w$, the
points $w_j$ have more preimages than $w$. For each $w_j$, call
one such preimage $p_j^* \not \in \cup_i V_i$. By properness, the
sequence $\{p_j^*\}$ must accumulate to a point $w_{\infty}$, and,
again by continuity, $w_{\infty}$ must be a preimage of $w$ which
does not belong to the interior of $\cup_i V_i$: this gives rise
to a new preimage of $w$, a contradiction. \qed

A proper excellent function $F: \RR^2 \to \RR^2$ is {\it nice}
if the following two conditions hold:
\begin{itemize}
\item{any point $y$ in $F(C)$ is the image of at most two critical points;}
\item{if $q$ is the image of two critical points $p_1$ and $p_2$
then both are folds and the tangent lines to $F(C)$ at $q$
corresponding to $p_1$ and $p_2$ are distinct.}
\end{itemize}
Points which are images of two critical points are {\it double
points}. Rather unsurprisingly, the generic excellent function is
nice, but we do not prove this technical result. From figure
\ref{fig:criticaldom}, the function $F_0$ is nice. Two distinct
tiles $A$ and $B$ for $F(C)$ are {\it adjacent} if their
boundaries share an arc of $F(C)$. In figure
\ref{fig:criticaldom}, tiles $Y_0$ and $Y_1$ are both adjacent to
$Y_2$ but not to each other.

\begin{lemma}
\label{lemma:neighbors}
Let $F: \RR^2 \to \RR^2$ be a nice function with critical set $C$;
the number of preimages of points in adjacent tiles for $F(C)$ differ by two.
\end{lemma}

\proof
Take $w$ to be the image of a fold $p$
belonging to a common boundary arc of adjacent tiles $A$ and $B$.
Let $p_1,\ldots,p_k$ be the preimages of $w$, with $p = p_1$.
As in the proof of the previous lemma,
we take disjoint open neighborhoods $V_i$ of $p_2, \ldots, p_k$
not containing $p$ which are taken homeomorphically by $F$
to an open neighborhood $W$ of $w$.
Take points $w_A \in A \cap W$, $w_B \in B \cap W$:
there will be $k-1$ preimages of $w_A$ and $w_B$ in the neighborhoods $V_i$.
Now, from the behavior of $F$ near $p$,
either $w_A$ or $w_B$ has two additional preimages close to $p$.
\qed

How can we obtain the {\it sense of folding}, i.e., on which of the 
two adjacent tiles for $F(C)$ do points have more preimages?
One way is to look at images of cusps:
from figure \ref{fig:cusp}, points inside the wedge
have more preimages than points outside it.

\section{Covering maps and the flower}

\begin{figure}[ht]
\vglue 11pt
\begin{center}
\epsfig{height=120mm,file=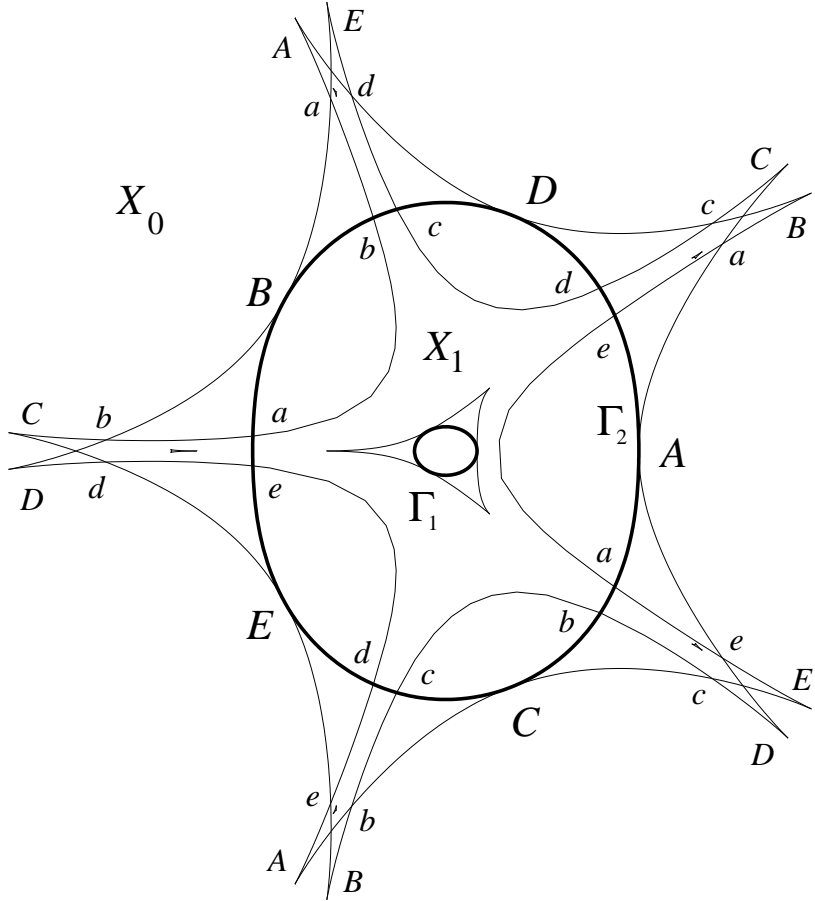}
\end{center}
\caption{\capsize The flower $F_0^{-1}(F_0(C))$}
\label{fig:flor}
\end{figure}

We now split the domain of a nice function $F$ in regions on which
$F$ behaves in a very simple fashion. More precisely, we consider
the tiles for $F^{-1}(F(C))$, the {\it flower} of $F$. Figure
\ref{fig:flor} shows the flower of $F_0$.

The two critical curves $\Gamma_1$ and $\Gamma_2$ are of course
part of the flower and are drawn thicker. The labels indicate
preimages of special points in the image, given in figure
\ref{fig:criticaldom}. The local behavior at the eight cusps of
$F_0$ is in agreement with figure \ref{fig:cusp}. Cusps in
$\Gamma_1$ are effective in the annulus between critical curves
and cusps in $\Gamma_2$ are effective in the region outside
$\Gamma_2$. Notice the five small preimages of $F_0(\Gamma_1)$ in
the five petal-like tiles for the flower: these are indeed
curvilinear triangles, as a zoom would show
(one is shown in figure \ref{fig:zoom}).

The tiles for the flower (resp. for $F(C)$) will be labelled $X_i$
(resp. $Y_j$). As we shall see, for many $X_i$, $F$ is a
diffeomorphism from $X_i$ to some $Y_j$, extending to a
homeomorphism between the closures $\overline{X_i}$ and
$\overline{Y_j}$. For the function $F_0$, for example, the only
exceptions are $X_0$ and $X_1$, indicated in figure
\ref{fig:flor}. It turns out that each point of $Y_0$
(see figure \ref{fig:criticaldom}) has $3$ preimages, all of them in $X_0$:
this is in agreement with lemma \ref{lemma:count} and
the fact that $F_0$ at infinity looks like $z \mapsto z^3$.
Furthermore, points in the boundary of $Y_0$ also have $3$ preimages
in the boundary of $X_0$ but may have other preimages elsewhere:
this can be checked by reading the labels
in figures \ref{fig:criticaldom} and \ref{fig:flor}.
Similarly, points in $\overline{Y_1}$ have $2$ preimages in
$\overline{X_1}$. Still, the restrictions $F: X_i \to Y_i$, $i = 0, 1$,
are examples of {\it covering maps},
a concept whose basic properties we now review
(\cite{Massey} and \cite{Munkres} are excellent references).

Take $X$ and $Y$ to be open nonempty, connected subsets of $\RR^2$:
in our examples, $X$ and $Y$ will be tiles $X_i$ and $Y_j$.
The continuous function $\Pi: X \to Y$ is a {\it covering map} if, for
any $y \in Y$, there exists an open neighborhood $V \subset \RR^2$
of $y$, such that $V \cap Y$ is connected and, for any connected
component $Z$ of $\Pi^{-1}(V \cap Y)$, the restriction $\Pi: Z \to V
\cap Y$ is a homeomorphism.

\begin{prop}
\label{prop:coverflower} Let $F: \RR^2 \to \RR^2$ be a nice
function with critical set $C$. Let $X_i$ and $Y_j$ be the tiles
for the flower $F^{-1}(F(C))$  and $F(C)$. Then the image of each
tile $X_i$ is a tile $Y_j$, the restriction $F: {X_i} \to {Y_j}$
is a covering map and $F: \overline{X_i} \to \overline{Y_j}$ is
locally injective.
\end{prop}

It is not always true that $F: \overline{X_i} \to \overline{Y_j}$
is injective: the boundary may contain two regular preimages of a
double point.

\proof Since $F(X_i) \subseteq \RR^2 - F(C)$ is connected, it is
contained in a single $Y_j$. Our proof of lemma \ref{lemma:count}
shows that the number $k$ of preimages under $F$ in $X_i$ is the
same for any point $y \in Y_j$ (and therefore $k > 0$). The
remaining argument is standard: given $y \in Y_j$, let $x_1,
\ldots, x_k \in X_i$ be its preimages (lemma \ref{lemma:finite}).
These are all regular points: by the inverse function theorem there are
disjoint open neighborhoods $U_1, \ldots, U_k \subset X_i$ of
$x_1, \ldots, x_k$ taken diffeomorphically to $V_1, \ldots, V_k$.
Take $V$ to be a small ball centered on $y$ contained in $V_1 \cap
\ldots \cap V_k$. This is the neighborhood of $y$ requested in the
definition of covering map. Indeed, since the number of preimages
is constant, there are no other preimages of $V$ outside
$U_1 \cup \ldots \cup U_k$. \qed

% We now apply these lemmas to our proper excellent function $F_0$.
% Suppose that we know that $\gamma_{\rho_1}$ and $\gamma_{\rho_2}$
% form the critical set $C$ of $F_0$.

% From lemma \ref{lemma:homeo} and the injectivity of $F_0$ on $\gamma_{\rho_1}$
% (which is checked numerically),
% we learn that $F_0$ is a homeomorphism between the disks
% $D_{\gamma_{\rho_1}}$ and $D_{\Gamma_{\rho_1}}$.
% Other components of $\RR^2 - C$ do not behave so nicely:
% instead, we consider components of $\RR^2 - F_0(C)$.

\begin{figure}[ht]
\vglue 11pt
\begin{center}
\epsfig{height=50mm,file=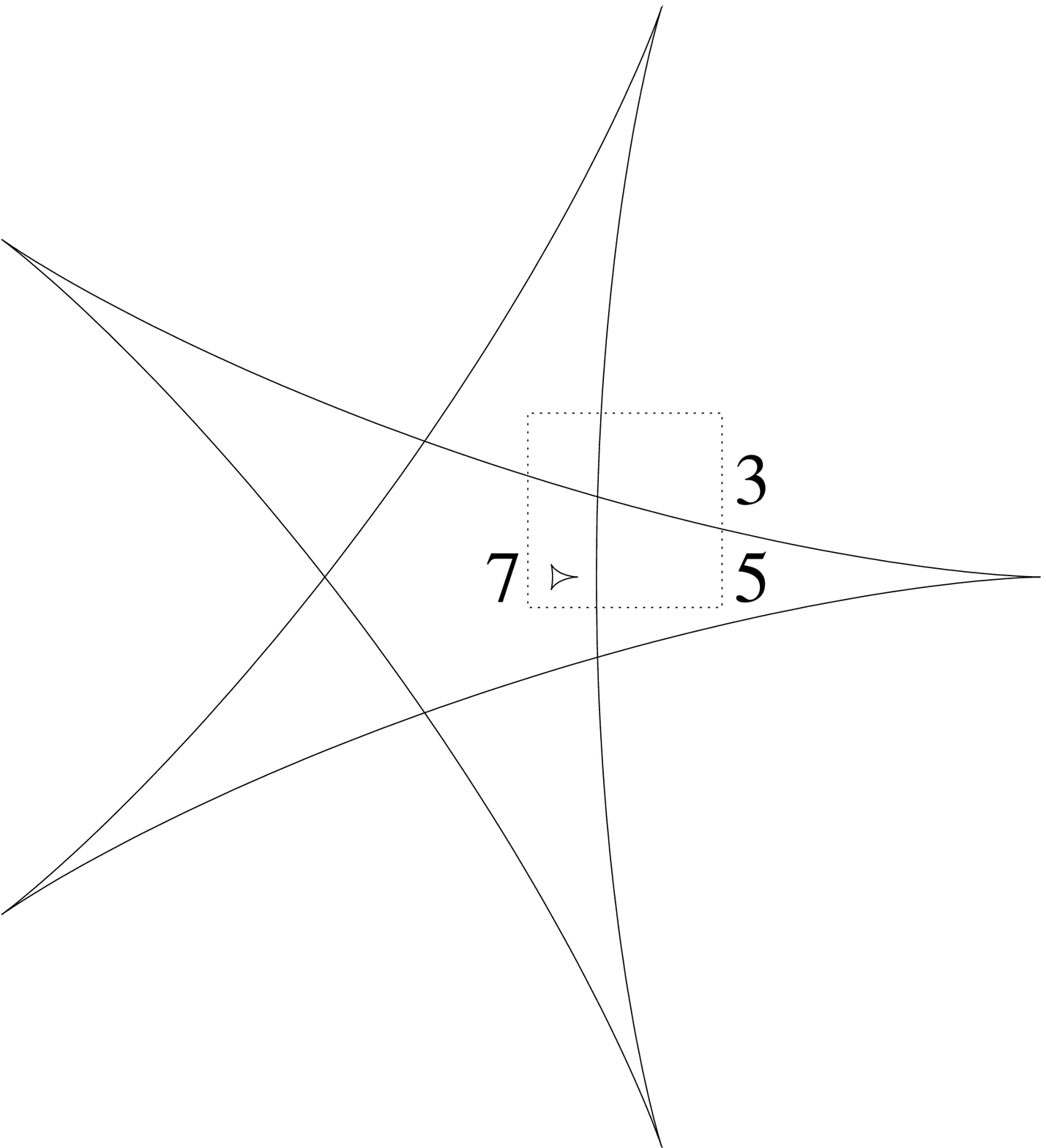}
\quad
\epsfig{height=50mm,file=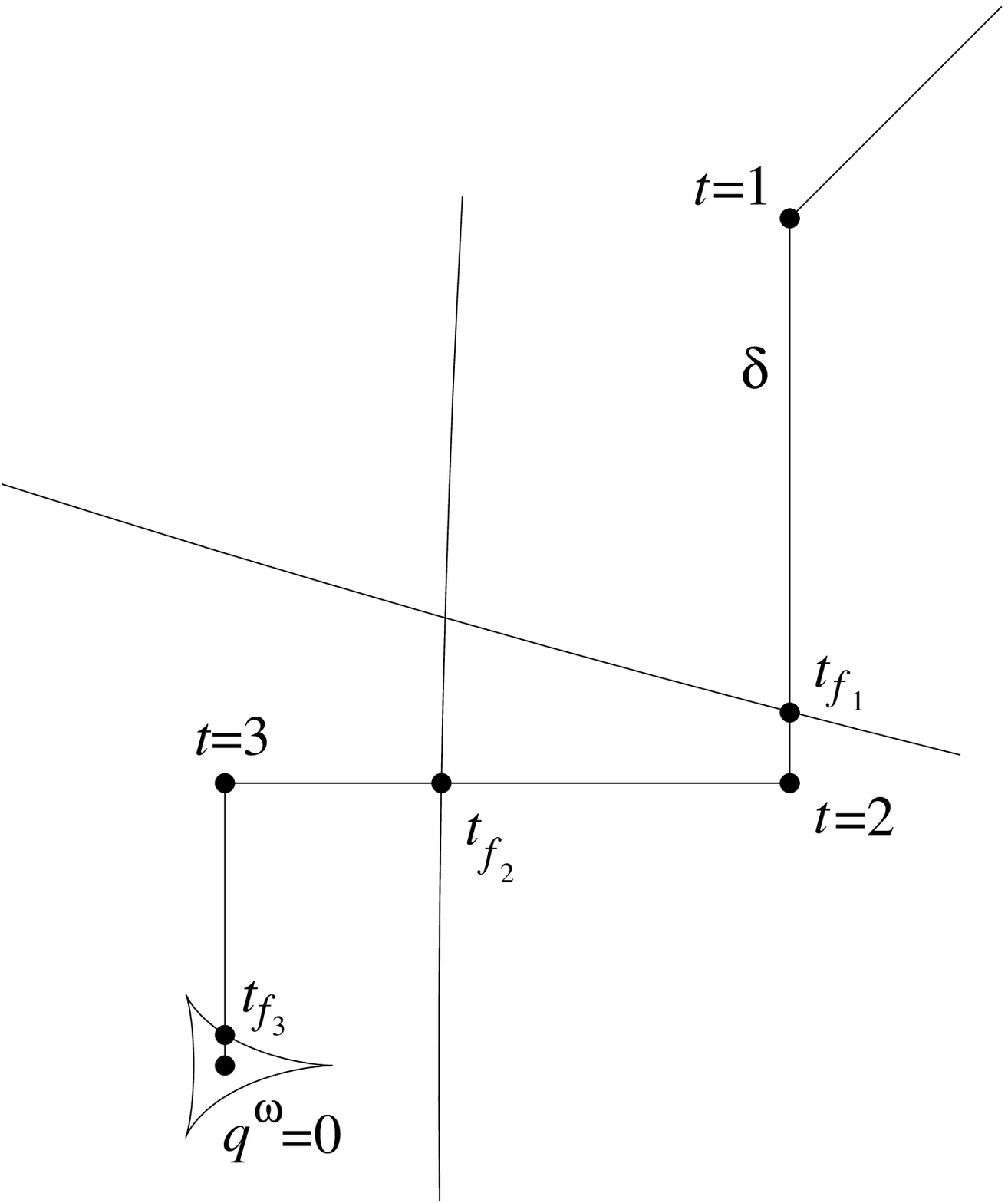}
\end{center}
\caption{\capsize Counting and computing preimages}
\label{fig:3to9}
\end{figure}

From the behavior of $F_0$ near infinity,
elements of large absolute value in the image
of $F_0$ have exactly three preimages.
Now, by lemmas \ref{lemma:count} and \ref{lemma:neighbors}
(using cusps to determine the sense of folding,
as suggested at the end of section 5),
we learn that the number of preimages in the tiles for
$F_0(C)$ vary as indicated in the left part of figure \ref{fig:3to9}.
The origin, which is at the very center of the innermost tile, has 9 preimages.
We can actually compute these preimages,
as we shall discuss in the next section:
they are the three conjugate pairs $1.864148 \pm 1.450656\;i$,
$-0.818866 \pm 2.665700\;i$, $0.204718 \pm 0.319589\;i$
and the three real numbers $0$, $-0.5$ and $-2$.

% This is accomplished by computing the inverse $F_0^{-1}(\zeta([0,4]))$
% under $F_0$ of the image of a path $\zeta: [0,4] \to \RR^2$
% joining $w_0 = \zeta(0)$ to $w_4 = 0 = \zeta(4)$,
% shown in the right part of figure \ref{fig:3to9}
% (the dotted square on the left bounds the region magnified on the right).
% More details on this inversion algorithm will be given in the next section.
% Figure \ref{fig:flowerpath}
% shows these preimages:
% the twelve endpoints of the six thicker curves
% are the three preimages of $w_0$ and the nine preimages of $w_4$.
% Proposition \ref{prop:coverflower} and the normal form of $F_0$ at folds
% guarantee that the connected components of $F_0^{-1}(\zeta([0,4]))$
% are homeomorphic to a closed interval.

We shall make use of the so called universal cover
of an open subset of the plane, as in following classical result.

\begin{theo}
\label{theo:universal}
Let $X \subset \RR^2$ be a nonempty connected open set;
then there exists a covering map $\Pi: \RR^2 \to X$.
\end{theo}

A sketch of proof could be as follows.
Consider $X \subset \CC$: if $X = \CC$ or $X = \CC - \{z_0\}$
set $\Pi(z) = z$ or $\Pi(z) = z_0 + \exp(z)$.
Otherwise, take $z_0 \in X$ and let $\Delta = \{ z \in \CC \; | \; |z| < 1 \}$:
clearly, $\Delta$ and $\RR^2$ are diffeomorphic.
Let $\cal F$ be the (nonempty) class of holomorphic functions
$f: \Delta \to X$, $f(0) = z_0$, $f'(0) > 0$.
In $\cal F$, there exists a function $f_0$ with maximum 
derivative at the origin.
Existence, uniqueness and the fact that $f_0 = \Pi$ is a covering map
follow as in the proof of the Riemann mapping theorem in
\cite{Ahlfors} or \cite{Rudin}.
With this proof, the theorem above is a special case
of the uniformization theorem (sections 3.2 and 3.3 of \cite{Lehto}).

\section{Computing preimages}

To compute preimages we use {\it continuation methods}, an example
of which we now describe (see \cite{AG} for more).
Let $F:\RR^2 \to \RR^2$ be a
nice function with critical set $C$, and take $p^\alpha$ to be a
regular point with image $q^\alpha = F(p^\alpha)$.
For a point $q$ sufficiently close to $q^\alpha$,
Newton's method computes the only preimage $p$ of $q$ near $p^\alpha$
by solving $F(p)=q$, taking $p^\alpha$ as the initial iteration.
Suppose now that we want to
compute a preimage of a point $q^\omega$ which is rather far from $q^\alpha$.
We draw a smooth parametrized arc $\delta: [0,1] \to \RR^2$
with $\delta(0) = q^\alpha$, $\delta(1) = q^\omega$ and
try to obtain points along a continuous path $\gamma: [0,1] \to \RR^2$
with $\gamma(0) = p^\alpha$, $F(\gamma(t)) = \delta(t)$.
More precisely, set $t_0 = 0 < t_1 < t_2 < \cdots < t_N = 1$
and try to compute $\gamma(t_{i+1})$ by solving
$F(\gamma(t_{i+1})) = \delta(t_{i+1})$ taking $\gamma(t_i)$
as initial condition for Newton's method.
If $\delta$ does not intersect $F(C)$ and the distances $t_{i+1} - t_i$
are taken to be sufficiently small then the method is guaranteed to
obtain $p^\omega = \gamma(1)$, a preimage of $q^\omega$.
This follows from the properness of $F$ combined with
the Newton-Kantorovich theorem (theorem 12.6.2, page 421, \cite{OR}).
If $\delta$ crosses $F(C)$, this continuation method may fail.
For instance, if $\delta(t_f) = F(p_f)$, where $p_f$ is a fold point,
as in figure \ref{fig:fold},
$\delta(t)$ belongs to the solid part of $\delta$ for $t < t_f$
(i.e., $\delta(t)$ belongs to the tile for $F(C)$ adjacent to $F(p_f)$
with the larger number of preimages)
and $\gamma(t)$ approaches $p_f$ when $t$ tends to $t_f$
then any continuation method ought to fail:
the dotted part of $\delta$ has no preimage near $p_f$
and no continuous function $\gamma$ with the required properties exists.

% We now consider the situation when $\delta$ is allowed to cross $F(C)$.
% If the intermediate points $q_i$ are chosen carefully,
% the process of obtaining $p_{i+1}$ from $p_i$ only breaks down
% if $p_i$ in near a critical point $p_f$:
% this only happens if $q_f = F(p_f)$
% is an intersection point of $\delta$ and $F(C)$.
% We assume that such intersections are generic.
% More precisely, $q_f$ is the image of a unique critical point $p_f$
% which is a fold; also, $\delta$ is transversal to $F(C)$ at $q_f$
% (i.e., their tangent vectors are linearly independent).
% The local form of $F$ at $p_f$, given in figure \ref{fig:fold},
% indicates what occurs in this case:
% the preimage of the arc $\delta$ hits the critical set $C$ of $F$ at $p_f$
% and there are no preimages of $\delta$ beyond $q_f$ near $p_f$.
% On the other hand, there is another preimage of points
% slightly before $q_f$ on the other side of $C$ near $p_f$.
% Thus, near $q_f$, the arc
% $\delta$ splits in two stretches $\delta_-$ and $\delta_+$ on
% opposite sides of $q_f$: $\delta_-$ has two preimages near $p_f$
% and $\delta_+$ has none.
% In particular, the inversion by continuation
% of the arc $\delta$ from $q_i$ to $q_{i+1}$
% fails to produce $p_{i+1}$ from $p_i$ precisely
% when $\delta$ crosses $F(C)$ from a region with more preimages
% to a region with fewer preimages.

We now consider the problem of computing all preimages of
a point $q^\omega \not\in F(C)$.
Assume that there exists $q^\alpha \not\in F(C)$
for which all preimages $p_1^\alpha, \ldots, p_n^\alpha$ are known.
Draw a piecewise smooth arc $\delta$ from $q^\alpha$ to $q^\omega$
which crosses $F(C)$ transversally at simple images of folds:
our strategy is to start with the set of all preimages of $q^\alpha$
and obtain all preimages of $q^\omega$
by an extension of a standard continuation method along $\delta$.
We may assume by induction that $\delta$ is smooth and intersects
$F(C)$ exactly once at $\delta(t_f) = q_f$.
Continuation along $\delta$ starting at each $p_i^\alpha$
tries to obtain paths $\gamma_i$ with $\gamma_i(0) = p_i^\alpha$,
$F(\gamma_i(t)) = \delta(t)$.
As we saw above, if $\delta$ crosses $F(C)$ from a tile
with more preimages to a tile with fewer preimages
then two of the paths $\gamma_i$ will collide at $p_f$
and will not be defined for $t > t_f$:
that is not a problem for us since the remaining paths
will still provide us with all the $n - 2$ preimages of $q^\omega$.
This scenario is reversed 
if $\delta$ crosses $F(C)$ from a tile with fewer preimages
(dotted in figure \ref{fig:fold}) to a tile with more preimages:
two new arcs are born at $p_f$.
More precisely, two distinct paths $\gamma_{n+1}$ and $\gamma_{n+2}$
from $[t_f,1]$ to $\RR^2$ exist with $\gamma_i(t_f) = p_f$
and $F(\gamma_i(t)) = \delta(t)$ for $i = n+1, n+2$, $t \ge t_f$.
These paths are quite removed from any of the $n$ preimages
$\gamma_i(t_f - \epsilon)$, $i = 1,\ldots,n$,
of $\delta(t_f - \epsilon)$ (for a small $\epsilon > 0$)
and could not possibly be obtained from these by a (local) continuation method.
Also, since the Jacobian $DF(p_f)$ is not invertible,
$p_f$ (which we know, since we previously obtained the critical curves)
is not acceptable as an initial condition for Newton's method
to solve $F(p) = q$ in $p$.
Instead, we compute a unit generator $v$ for $\ker DF(p_f)$
and set $p_{n+1} = p_f + sv$, $p_{n+2} = p_f - sv$
(for a small positive real number $s$)
and $q_i = F(p_i)$ ($i = n+1, n+2$).
From the normal form,
each $q_i$ is now not too far from $\delta(t_f + \epsilon)$
(for some small $\epsilon$) and can be connected
to it by an auxiliary arc $\delta_i$ which does not intersect $F(C)$:
our continuation method now obtains $\gamma_i(t_f + \epsilon)$
by following $\delta_i$, starting with $p_i$.
Recall that the preimage of any smooth curve $\delta$
crossing $F(C)$ transversally at $q_f$ is tangent to $v$ at $p_f$
(figure \ref{fig:fold}).

\begin{figure}[t]
\vglue 11pt
\begin{center}
\epsfig{height=120mm,file=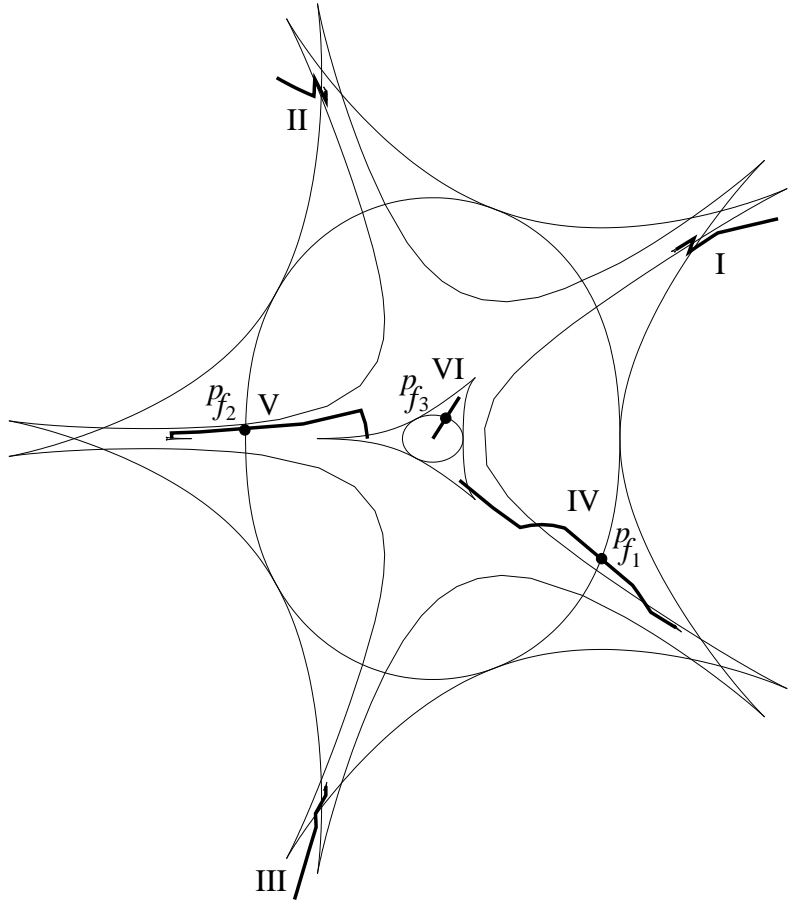}
\end{center}
\caption{\capsize Inverting a path}
\label{fig:flowerpath}
\end{figure}

Let us now go back to our basic example and see how our program
obtains the nine preimages of $0$ under $F_0$.
First it computes the critical set $C$ of $F_0$ and its image,
presented in figure \ref{fig:criticaldom}.
Next, it obtains the three preimages of
a remote point $q^\alpha$ (see figure \ref{fig:3to9}).
This is rather simple: for complex numbers $z$ of large absolute value
the function $F_0(z)$ is similar to $z \mapsto z^3$
and the three complex cube roots of $q^\alpha$
are good initial conditions for a Newton-like method to solve
$F_0(p^\alpha_i) = q^\alpha$, $i = 1,2,3$.
The three preimages lie in the regions indicated by
the Roman numerals I, II and III in figure \ref{fig:flowerpath}.

A path $\delta: [0,4] \to \RR^2$
from $q^\alpha = \delta(0)$ to $0 = \delta(4)$
% (as in figure 6.1)
(as in figure \ref{fig:3to9}) was constructed
as a juxtaposition of four smooth paths
defined on intervals with integer endpoints.
% of the kind discussed in the previous paragraphs:
The first ($0 \le t \le 1$) is the only one that does not cross $F_0(C)$.
Notice that the number of preimages is increasing along this path.

\begin{figure}[ht]
\vglue 11pt
\begin{center}
\epsfig{height=60mm,file=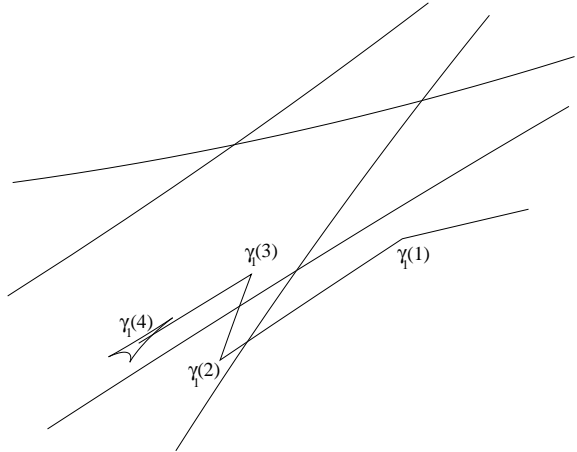}
\end{center}
\caption{\capsize Zoom on region I of figure 8}
\label{fig:zoom}
%\ref{fig:flowerpath}}
\end{figure}

Three paths $\gamma_i: [0,4] \to \RR^2$, $i = 1,2,3$,
were obtained by a continuation method starting from
$\gamma_i(0) = p^\alpha_i$.
The path $\gamma_1$, which is in region I,
is presented in figure \ref{fig:zoom};
$\gamma_2$ and $\gamma_3$ are in regions II and III
in figure \ref{fig:flowerpath}.
The whole inversion procedure from $\gamma_i(0)$ to $\gamma_i(4)$
does not cross a critical curve of $F_0$,
and three solutions to the equation $F_0(z) = 0$ are obtained:
$\gamma_1(4) \approx (1.864148, 1.450656)$,
$\gamma_2(4) \approx (-0.818866, 2.665700)$ and
$\gamma_3(4) \approx (-0.818866, -2.665700)$.

The paths obtained by continuation within regions I, II and III do not
notice anything unusual at $t_{f_1}$,
the first intersection between $\delta$ and $F_0(C)$.
As we saw, however, two new arcs $\gamma_4$ and $\gamma_5$ are born
at the critical point $p_{f_1}$ for which $F_0(p) = \delta(t_{f_1})$.
The program identified $p_{f_1} = \gamma_4(t_{f_1}) = \gamma_5(t_{f_1})$,
which turns out to lie in the outer critical curve,
and obtained by continuation from $p_{f_1}$ two new paths, lying in region IV.
Similarly, two new paths are born at $t_{f_2}$ from the fold $p_{f_2}$
(they are in region V) and yet two more at $t_{f_3}$ from $p_{f_3}$
(in region VI).

These computations rely heavily on the assumption that
the critical set of $F_0$ has been correctly identified.
This is the same issue raised at the end of section 4;
we next introduce topological tools to tackle this problem.

\section{Rotation numbers}

We remind the reader of a few facts concerning winding numbers
(for a more complete exposition, see \cite{CS}, sections 17 to 27).
For a continuous function $\phi: [a,b] \to \RR^2$ and $p \in
\RR^2$, $p$ not in the image of $\phi$, define a continuous
{\it argument function} $\theta_p: [a,b] \to \RR$ such that
\[ \phi(t) = |\phi(t) - p| (\cos\theta_p(t), \sin\theta_p(t)) \]
for all $t \in [a,b]$. The argument function is unique up to an
additive constant of the form $2\pi n$ and the {\it angle swept}
by $\phi$ with respect to $p$, $A_p(\phi,p) = \theta_p(b) -
\theta_p(a)$, is well defined. Parametrize the standard unit
circle by $e: [0,2\pi] \to \Ss^1$, where $e(t) = (\cos t, \sin t)$.
For a closed curve $c: \Ss^1 \to \RR^2$ with $p$ not in the image
of $c$ we have that $W(c,p) = A_p(c \circ e,p)/(2\pi)$ is an
integer which we call the {\it winding number} of $c$ around $p$.

A {\it homotopy} between continuous functions $\phi_0: [a,b] \to
\RR^2$ and $\phi_1: [a,b] \to \RR^2$ is a continuous function
$\Phi: [0,1] \times [a,b] \to \RR^2$ with $\Phi(0,t) = \phi_0(t)$
and $\Phi(1,t) = \phi_1(t)$.
The winding number around $p$ is invariant under homotopy
provided all curves are closed and avoid the point $p$.
More precisely, let $\Phi: [0,1] \times [a,b] \to \RR^2$ be a homotopy
between $\phi_0: [a,b] \to \RR^2$ and $\phi_1: [a,b] \to \RR^2$,
so that $\Phi(s,a) = \Phi(s,b)$ for all $s$, for which $p$ is not
in the image of $\Phi$. Then we must have $W(\phi_0,p) =
W(\phi_1,p)$ (theorem 25.1 in \cite{CS}).

Of special interest will be {\it rotation numbers}:
we present the basic results following \cite{H} and \cite{W1}.
A {\it parametrized regular closed curve} (in short, prc-curve)
is a $C^1$ function $c: \Ss^1 \to \RR^2$
with $(c \circ e)'(t) \ne 0$ for all $t \in [0,2\pi]$:
the image of $c$ is an oriented curve $\gamma$,
possibly with self-intersections.
The {\it rotation number} $r(c)$
is the winding number of $c'$ around $0$: $r(c) = W(c',0)$.

Two prc-curves $c_0$ and $c_1$ are {\it equivalent} if their
images as oriented curves are equal, or, more precisely,
if there exists an orientation preserving $C^1$ diffeomorphism
$\eta: \Ss^1 \to \Ss^1$ with $c_0 = c_1 \circ \eta$.
Equivalent curves have the same rotation number (\cite{W1}):
this allows the computation of the rotation number of a prc-curve
from the drawing of its image.

\begin{figure}[ht]
\vglue 11pt
\begin{center}
\epsfig{height=25mm,file=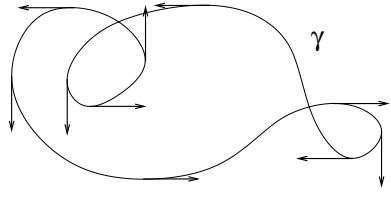}
\end{center}
\caption{\capsize Computing rotation numbers}
\label{fig:howtorotate}
\end{figure}

A recipe to compute $r(c) = r(\gamma)$ is the following.
Draw all horizontal and vertical tangent vectors to
the oriented curve $\gamma$ as in figure \ref{fig:howtorotate},
measure the oriented angles between neighboring vectors
(always equal to $0$, $\pi/2$ or $-\pi/2$)
add them all up and divide by $2\pi$.
In the figure, $r(\gamma) = 1$.
As another example, denote by $e_\rho$ a counterclockwise
parametrization of the circle of radius $\rho$ around the origin.
The curves in figure \ref{fig:3regs} are thus parametrized
by $F_0 \circ e_\rho$ for various values of $\rho$.
Their rotation numbers are $1$, $-2$ and $3$, respectively.

The following result, known as the {\it Umlaufsatz},
is somewhat harder to prove.

\begin{theo}[Hopf, \cite{H}]
\label{theo:rotatesimple}
If $c$ is an injective prc-curve then $r(c) = \pm 1$.
\end{theo}

\proof
Let $t_0 \in \Ss^1$ be a point maximizing $|c(t)|, t \in \Ss^1$.
Let $\ell$ be the line tangent to $\gamma$, the image of $c$, through $c(t_0)$.
By construction, $c(t_0)$ is the only common point
between $\ell$ and $\gamma$ and $\gamma$ is a subset
of the half-plane defined by $\ell$ containing the origin.
Without loss, set $t_0 = 0$, $\ell$
to be the horizontal line $y = -1$ and $(c \circ e)'(0) = (p,0)$
where $p > 0$: in this case, we show that $r(c) = 1$.
Let $T = \{ (s,t); 0 \le s \le t \le 2\pi \}$ and
$\psi: T \to \RR^2$ be defined by
\[ \psi(s,t) =
\begin{cases}
(c \circ e)'(t) & \textrm{for }s = t, \\
-(c \circ e)'(0) & \textrm{for }s = 2\pi \textrm{ and } t = 0, \\
\frac{\textstyle (c \circ e)(t) - (c \circ e)(s)}{\textstyle \min \{t - s, 1 + s - t\} } &
\textrm{otherwise.}
\end{cases} \]
Since $c$ is a $C^1$ function, $\psi$ is continuous;
injectivity of $c$ implies that $\psi$ is never zero.
The path $\phi_1(t) = \psi(0,t)$, $t \in [0,2\pi]$,
satisfies $\phi_1(0) = (p,0)$, $\phi_1(2\pi) = (-p,0)$
and $\phi_1(t)$ stays above the horizontal axis,
whence $\phi_1$ sweeps half a turn: $A(\phi_1,0) = \pi$.
Similarly, for $\phi_2(s) = \psi(s,2\pi)$, $s \in [0,2\pi]$,
we have $A(\phi_2,0) = \pi$.
Juxtapose $\phi_1$ and $\phi_2$ to define
$\phi_{12}: [0,2\pi] \to \RR^2$ where
$\phi_{12}(t) = \phi_1(2t)$ for $t \in [0,\pi]$
and $\phi_{12}(t) = \phi_2(2t - 2\pi)$ for $t \in [\pi,2\pi]$:
clearly, $A(\phi_{12},0) = 2\pi$.
The function $\psi$ can be viewed as a
homotopy between $(c \circ e)'$
(the restriction of $\psi$ to $\{(t,t),t \in [0,2\pi]\}$)
and $\phi_{12}$
(the restriction to
$\{(0,t), t \in [0,2\pi]\} \cup \{(s,2\pi), s \in [0,2\pi]\}$),
showing that $A((c \circ e)',0) = A(\phi_{12},0) = 2\pi$
and thus $r(c) = 1$. \qed

Two prc-curves $c_0$ and $c_1$ can be {\it deformed} into each
other if there exists a continuous function
$H: [0,1] \times [0,2\pi] \to \RR^2$
such that $H(s,t) = c_s(e(t))$ for all $s \in \{0,1\}$, $t \in [0,2\pi]$,
$\frac{\partial H}{\partial t}$ is continuous and nonzero
in $[0,1] \times [0,2\pi]$ and
$\frac{\partial H}{\partial t}(s,0) = \frac{\partial H}{\partial t}(s,2\pi)$
for all $s \in [0,1]$.
As the next theorem shows, this
is the appropriate concept of deformation on prc-curves, if we want to
preserve rotation number. For instance, the curves in figure
\ref{fig:3regs} do not admit deformations joining them.

\begin{theo}[Graustein and Whitney, \cite{W1}]
\label{theo:deform}
Two prc-curves $c_0$ and $c_1$
can be deformed into each other if and only if $r(c_0) = r(c_1)$.
\end{theo}

\proof The invariance of rotation number under deformation is a
corollary of the invariance of winding number under homotopy:
this proves one implication.

% More precisely, if $c_0, c_1: [0, 2\pi] \to \RR^2$
% are two $C^1$ closed curves for which there exists a homotopy
% $\Phi: [0,2\pi] \times [0,1] \to \RR^2 - 0$ with
% $\Phi(t,0) = c'_0(t)$, $\Phi(t,1) = c'_1(t)$ for all $t$,
% $\Phi(0,s) = \Phi(2\pi,s)$ for all $s$
% then $r(c_0) = r(c_1)$.
% In particular, rotation numbers do not depend on parametrization.
% This was of course suggested by the notation and is also the reason
% why the recipe above depends on the curve only.

Now, let $c_0$ and $c_1$ be prc-curves with $r(c_0) = r(c_1) = n$.
Reparametrize by arc length and change scale so that $|(c_i \circ e)'(t)| = 1$
for all $t \in [0,2\pi]$, $i = 0, 1$.
Let $X: [0,1] \times [0,2\pi] \to \Ss^1 \subset \RR^2$
be a continuous function with $X(i,t) = (c_i \circ e)'(t)$,
$X(s,0) = X(s,2\pi)$ and such that, for any $s \in [0,1]$,
the function $t \mapsto X(s,t)$ is not constant:
the existence of such $X$ is a standard topological fact,
but for completeness we provide an explicit construction.
Let $\theta_i: [0,2\pi] \to \RR$ be argument functions for $(c_i \circ e)'$:
we have $\theta_i(2\pi) - \theta_i(0) = 2\pi n$ for both values of $i$.
For $s \in [0,1]$, consider the segment joining $\theta_0$ and $\theta_1$:
$\tilde\theta_s(t) = (1-s) \theta_0(t) + s \theta_1(t)$.
If $n \ne 0$, $\tilde\theta_s$ is clearly not a constant function
and we take $X(s,t) = e(\tilde\theta_s(t))$.
For $n = 0$, take $\theta_{1/2}: [0,2\pi] \to \RR$ to be
an arbitrary continuous function with $\theta_{1/2}(0) = \theta_{1/2}(2\pi)$
which is not contained is the linear subspace generated by
$\theta_0$, $\theta_1$ and the constant function $1$.
Define $\theta_s$, $s \in [0,1]$,
by juxtaposing segments from $\theta_0$ to $\theta_{1/2}$
and from there to $\theta_1$;
take $X(s,t) = e(\theta_s(t))$.
Let \[ m(s) = \frac{1}{2\pi} \int_0^{2\pi} X(s,t)\,dt
\quad \textrm{and} \quad
Y(s,t) = X(s,t) - m(s). \]
Notice that $|m(s)| < 1$ and therefore $Y(s,t) \ne 0$ for all $s$ and $t$.
Also, for any given $s$, the integral of $Y(s,t)$ is $0$
so that $c_s$ defined by
\[ (c_s \circ e)(t) = \int_0^t Y(s,\tau)\,d\tau \]
is a prc-curve.
This is the required deformation.
\qed

% Rotation numbers are strongly related to the critical set.
% Indeed, in the previous example, the rotation number changes
% exactly when the circles $\eta_\rho$ trespass the critical curves.

Two different parametrizations of the same oriented
smooth curve $\gamma$, yielding two prc-curves,
have the same rotation and can therefore be deformed into each other.
We may therefore ask, without ambiguity,
whether two smooth curves $\gamma_0$ and $\gamma_1$
can be deformed into each other
(within the class of prc-curves $c_s: \Ss^1 \to \RR^2$):
this happens if and only if $r(\gamma_0) = r(\gamma_1)$.

\goodbreak

\section{Compatibility checks}

Suppose that we have detected some critical curves,
forming a certain subset $C_1$ of
the critical set $C$ of an excellent function $F$.
The propositions in this section provide global compatibility checks on $C_1$,
i.e., necessary (but not sufficient) conditions for $C_1 = C$.
We start with a technical lemma.

\begin{lemma}
\label{lemma:whitney2}
Let $U \subseteq \RR^2$ be a connected open set
and let $\gamma_0, \gamma_1$ be positively oriented
smooth simple closed curves bounding closed topological disks
$\Delta_0, \Delta_1 \subset U$.
Then $\gamma_0$ can be deformed to $\gamma_1$ within $U$,
i.e., the image of the deformation is contained in $U$.
\end{lemma}

\proof
Let $\Pi: \RR^2 \to U$ be a covering map
(theorem \ref{theo:universal}) and
consider $\Pi^{-1}(\Delta_0)$.
This set is a disjoint union of closed disks:
let $\tilde\Delta_0$ be one of them
and $\tilde\gamma_0$ its boundary.
Construct $\tilde\gamma_1$ similarly.
Since $\tilde\gamma_s$, $s = 0, 1$, are both simple curves,
$r(\tilde\gamma_0) = r(\tilde\gamma_1) = 1$ (theorem \ref{theo:rotatesimple})
and therefore $\tilde\gamma_0$ may be deformed to $\tilde\gamma_1$
(theorem \ref{theo:deform}).
Composing this deformation with $\Pi$ yields a deformation
from $\gamma_0$ to $\gamma_1$ contained in $U$, as desired.
\qed

\begin{prop}
\label{prop:rotatedisk}
Let $F: \RR^2 \to \RR^2$ be an excellent smooth function
with critical set $C$.
Let $\gamma$ be a positively oriented smooth simple closed curve
bounding a closed topological disk $\Delta$ with $\Delta \cap C = \emptyset$.
Then $F(\gamma)$, the image of $\gamma$ under $F$, is a smooth curve and
$r(F(\gamma)) = \sgn \det DF(p)$ for any $p \in \Delta$.
\end{prop}

Here $\sgn(x)$ is the usual sign function,
$\sgn(x) = 1$ (resp. $-1$) for $x > 0$ (resp. $x < 0$).

\proof
First notice that the result holds if $F$ is affine and the
curve $\gamma$ is a circle: $F(\gamma)$ is an ellipse and its
orientation is given by $\sgn \det DF(p)$.
Next consider an arbitrary $F$ and $p \not\in C$.
The affine map $\tilde F(v) = F(p) + DF(p) \cdot (v-p)$
is a $C^1$ approximation of $F$ around $p$:
thus there exists $\rho_0$ such that,
if $\gamma$ is a positively oriented circle of radius $\rho$ around $p$,
$0 < \rho < \rho_0$, then $|r(F(\gamma)) - r(\tilde F(\gamma))| < 1$
(the arguments of the tangent vectors are arbitrarily close
for small $\rho_0$).
Since rotation numbers are integers, $r(F(\gamma)) = r(\tilde F(\gamma))$.
Let now $\gamma_0$ be arbitrary
and $\gamma_1$ be a small round circle around some $p \in \Delta$:
use lemma \ref{lemma:whitney2} to deform $\gamma_0$ to $\gamma_1$
within the connected component of $\RR^2 - C$ containing $p$.
Compose this deformation with $F$ to conclude
that $r(F(\gamma_0)) = r(F(\gamma_1))$.
\qed

\begin{prop}
\label{prop:rotatepolydisk}
Let $F: \RR^2 \to \RR^2$ be an excellent smooth function
with critical set $C$.
Let $\gamma_0, \gamma_1, \ldots, \gamma_n$ be  positively oriented
smooth simple closed curves
bounding $\Delta_n$, a closed topological disk with $n$ holes
with $\Delta_n \cap C = \emptyset$.
Assume $\gamma_0$ to be the outer connected component
of the boundary of $\Delta_n$.
Let $s = \sgn\det DF(p)$, $p \in \Delta_n$.
Then
\[ r(F(\gamma_0)) = r(F(\gamma_1)) + \cdots + r(F(\gamma_n)) - s(n-1). \]
\end{prop}

In particular, if $n = 1$, we learn that if $r(F(\gamma_0)) = r(F(\gamma_1))$,
in agreement with lemma \ref{lemma:whitney2}.

\begin{figure}[ht]
\vglue 11pt
\begin{center}
\epsfig{height=25mm,file=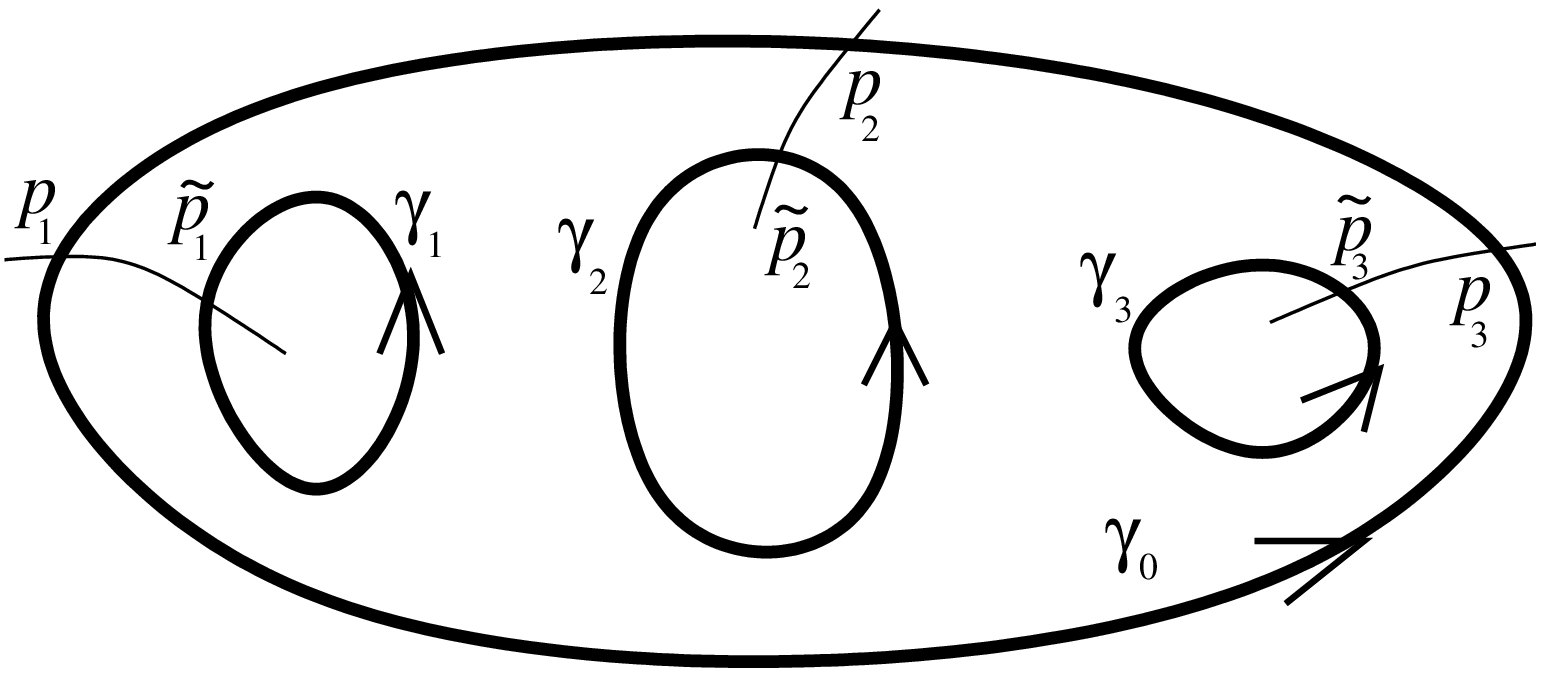}
\qquad
\epsfig{height=25mm,file=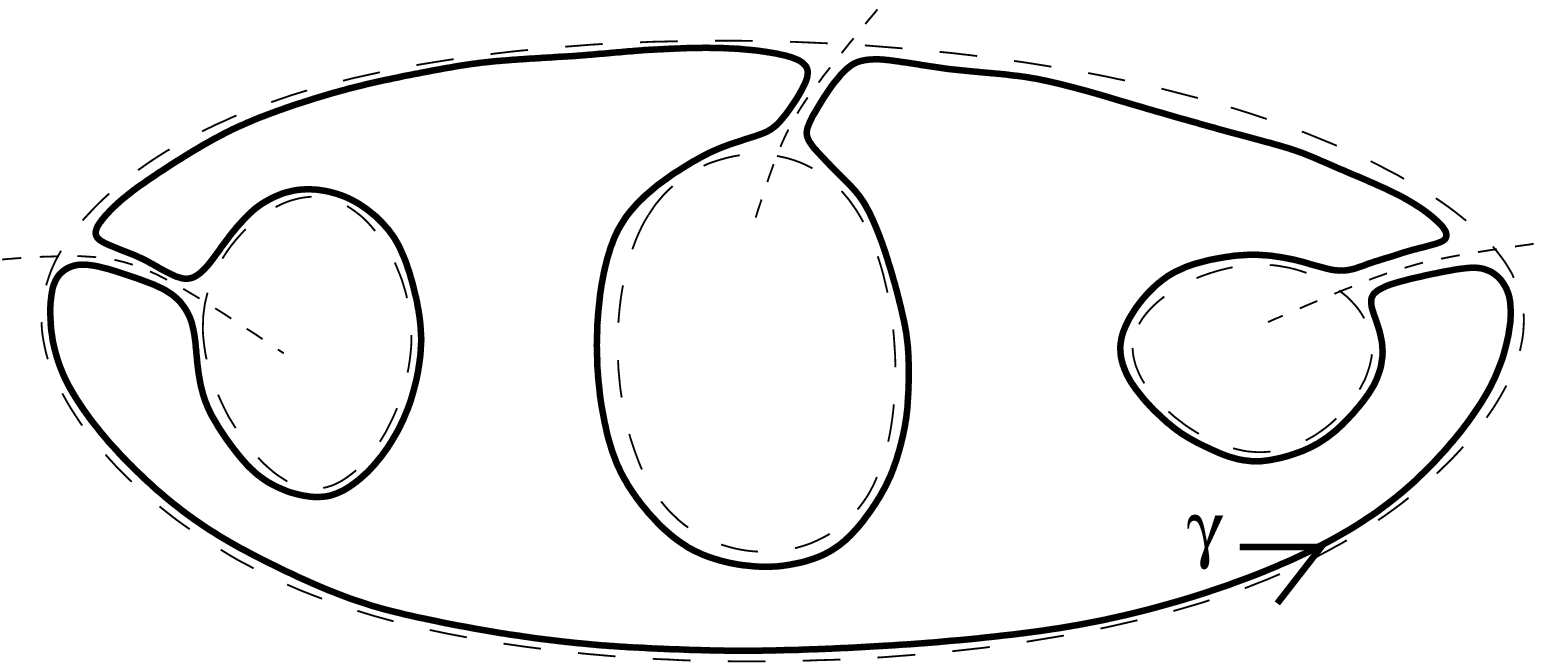}
\end{center}
\caption{\capsize Adding rotation numbers}
\label{fig:rotatepolydisk}
\end{figure}

\proof
Construct $n$ smooth disjoint arcs $\delta_1, \ldots, \delta_n$
such that $\delta_j$ crosses $\gamma_0$ and $\gamma_j$ transversally
at points $p_j$ and $\tilde p_j$, respectively,
as in figure \ref{fig:rotatepolydisk} (a).
Construct a simple prc-curve $\gamma$
in the interior of $\Delta_n - \cup_i \delta_i$
close to its boundary,
as indicated in figure \ref{fig:rotatepolydisk} (b).
% starting at an arbitrary point of $\gamma_0$,
% following it almost to the first $p_j$,
% introduce a smooth arc $\eta_{j,1}$ to then follow $\delta_j$
% almost to $\tilde p_j$,
% introduce another smooth arc $\eta_{j,2}$ to follow $\gamma_j$ backwards
% almost back to $\tilde p_j$, then another arc $\eta_{j,3}$,
% then $\delta_j$ almost back to $p_j$, then $\eta_{j,4}$,
% then $\gamma_0$ to the next point $p_{j'}$ and so on
% as in figure \ref{fig:rotatepolydisk} (b).
From proposition \ref{prop:rotatedisk}, $r(F(\gamma)) = s$.

On the other hand,
\[ r(F(\gamma)) =
r(F(\gamma_0)) - r(F(\gamma_1)) - \cdots - r(F(\gamma_n)) + sn. \]
Indeed, the parts of $\gamma$ close to some $\gamma_i$
contribute to $r(F(\gamma))$, up to a small error, with
$r(F(\gamma_0)) - r(F(\gamma_1)) - \cdots - r(F(\gamma_n))$.
Similarly, the parts of $\gamma$ near some $\delta_j$
contribute, again up to a small error,
with $0$ since arcs on either side of $\delta_j$
essentially cancel their contributions.
For each of the $2n$ intersections between some $\gamma_i$ and some $\delta_j$,
there are two small arcs of $\gamma$ which together contribute
with half a turn, more precisely, with approximately $s/2$.
Finally, the small errors cancel each other since both right and left
hand side are integers.
\qed

A cusp $q$ in a closed critical curve $\Gamma$ is
an {\it inner} (resp. {\it outer}) cusp
if it is effective on the bounded (resp. unbounded) component
of $\RR^2 - \Gamma$.

\begin{prop}
\label{prop:rotatecritical}
Let $F: \RR^2 \to \RR^2$ be an excellent smooth function with critical set $C$.
Let $A$ be a closed annulus
containing a single critical curve $\Gamma = A \cap C$;
set $\gamma_\textit{in}$ and $\gamma_\textit{out}$
to be the positively oriented simple closed 
components of the boundary of $A$, assumed to be smooth.
Let $k_\textit{in}$ and $k_\textit{out}$
be the number of inner and outer effective cusps
on $\Gamma$ and let $s_\textit{in} = \sgn\det DF(p_\textit{in})$,
$p_\textit{in} \in \gamma_\textit{in}$
and $s_\textit{out} = -s_\textit{in} = \sgn\det DF(p_\textit{out})$,
$p_\textit{out} \in \gamma_\textit{out}$.
Then
\[ r(F(\gamma_\textit{out})) =
r(F(\gamma_\textit{in})) + s_\textit{in} k_\textit{in} +
s_\textit{out} k_\textit{out}. \]
\end{prop}

\proof
Given proposition \ref{prop:rotatepolydisk} (with $n = 1$) we may assume
$\gamma_\textit{in}$ and $\gamma_\textit{out}$ to be very near $\Gamma$
and for their tangent vectors to be likewise near
the tangent vectors to $\Gamma$.
We first deform $\gamma_\textit{in} = \gamma_0$ into $\gamma_1$,
a curve which coincides with $\gamma_\textit{out}$ except
in small neighborhoods of cusps.
If in the process the tangent vectors to intermediate curves $\gamma_s$,
$s \in [0,1]$, are kept almost parallel to the tangent vectors to $\Gamma$,
they will never lie in the kernel of $DF$
(which, due to the normal form at folds, is never tangent to $\Gamma$)
and thus $F(\gamma_s)$ are all regular and
$r(F(\gamma_1)) = r(F(\gamma_\textit{in}))$.

Let $\gamma_2$ be a curve which coincides with $\gamma_1$ everywhere
except in the neighborhood of a cusp $p_c$, where $\gamma_2$ coincides
with $\gamma_\textit{out}$.
We now compare $r(F(\gamma_2))$ and $r(F(\gamma_1))$.
The region where $\gamma_1$ and $\gamma_2$ do not coincide lies
in a small neighborhood of $p_c$ and we may therefore use the normal form
at cusps: 
$r(F(\gamma_2)) - r(F(\gamma_1)) = \sgn\det DF(p)$,
where $p$ is in the region where the cusp $p_c$ is effective.
Repeating this process for the other curves yields the desired result.
\qed

As an application, we perform some tests on $F_0$. Suppose that we know that
the critical set $C$ contains (at least) the curves $\Gamma_1$ and $\Gamma_2$
as indicated in figure \ref{fig:criticaldom}.
Let $D_{\Gamma_1}$ and $D_{\Gamma_2}$ be the topological
disks bounded by these curves.
Consider four simple, positively oriented closed
curves $\gamma_{\textit{in},1}$, $\gamma_{\textit{out},1}$,
$\gamma_{\textit{in},2}$, $\gamma_{\textit{out},2}$,
on both sides of the critical curves,
as in proposition \ref{prop:rotatecritical}.
From the knowledge of the images of these four curves we learn that
$r(F_0(\gamma_{\textit{in},1})) = 1$,
$r(F_0(\gamma_{\textit{out},1})) = -2$,
$r(F_0(\gamma_{\textit{in},2})) =  -2$ and
$r(F_0(\gamma_{\textit{out},1})= 3$.
Cusps on $\Gamma_j, j=1,2$ are effective on the outside of $D_{\Gamma_j}$;
these values are therefore in agreement with
proposition \ref{prop:rotatecritical}.

If the rotation number of $F_0(\gamma_{\textit{in},1})$ were different from $1$,
we would learn from proposition \ref{prop:rotatedisk} that there had to be
additional critical curves in $D_{\Gamma_1}$.
Similarly, proposition \ref{prop:rotatepolydisk} would indicate the presence
of critical curves in the annulus bounded by $\Gamma_1$ and $\Gamma_2$
if the rotation numbers $r(F_0(\gamma_{\textit{out},1}))$ and
$r(F_0(\gamma_{\textit{in},2}))$ were different.
Finally, again by proposition \ref{prop:rotatepolydisk},
$r(F_0(\gamma_{\textit{out},2}))= 3$
is compatible with the behavior of $F_0$ at infinity.

Given a nice function $F$ and a subset $C_1$ of the critical set $C$,
our criteria never guarantee that $C_1 = C$.
There is a more complicated theorem which provides
necessary and sufficient conditions for the existence
of a nice function $F_1$ coinciding with $F$ in a neighborhood
of $C_1$ and having $C_1$ as critical set.
Theorem 3.1 in \cite{MST1}
(or theorem 1.6 in \cite{MST2} for the simpler case of bounded critical sets)
makes use of the ingredients used in this text,
combined with an additional tool
from combinatorial topology---Blank-Troyer theory (\cite{Poe}, \cite{Troyer}).
Whether $C_1$ is the critical set of $F$ is something that
cannot be resolved without invoking completely different methods.
The difficulty is already evident in one dimension: how do you know that
your favorite numerical method has found all the roots of, say,
the real function $f(x) = x$?
If one is only entitled to a finite number of evaluations of a function
and its derivatives,
one will never know what happens on very small scales or
at very remote points.

\section{Other examples}

\begin{figure}[ht]
\vglue 11pt
\begin{center}
\epsfig{height=40mm,file=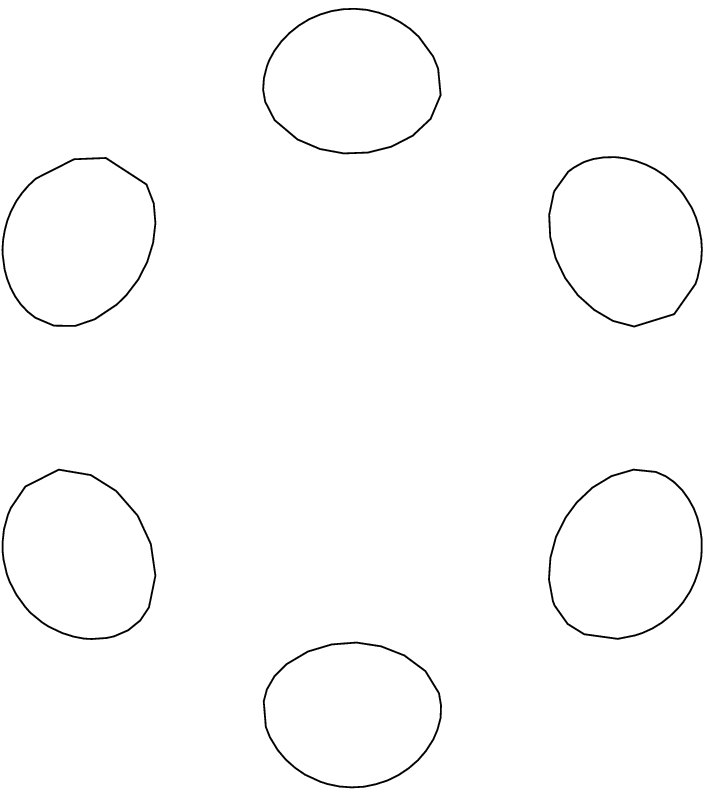}
\qquad \qquad
\epsfig{height=40mm,file=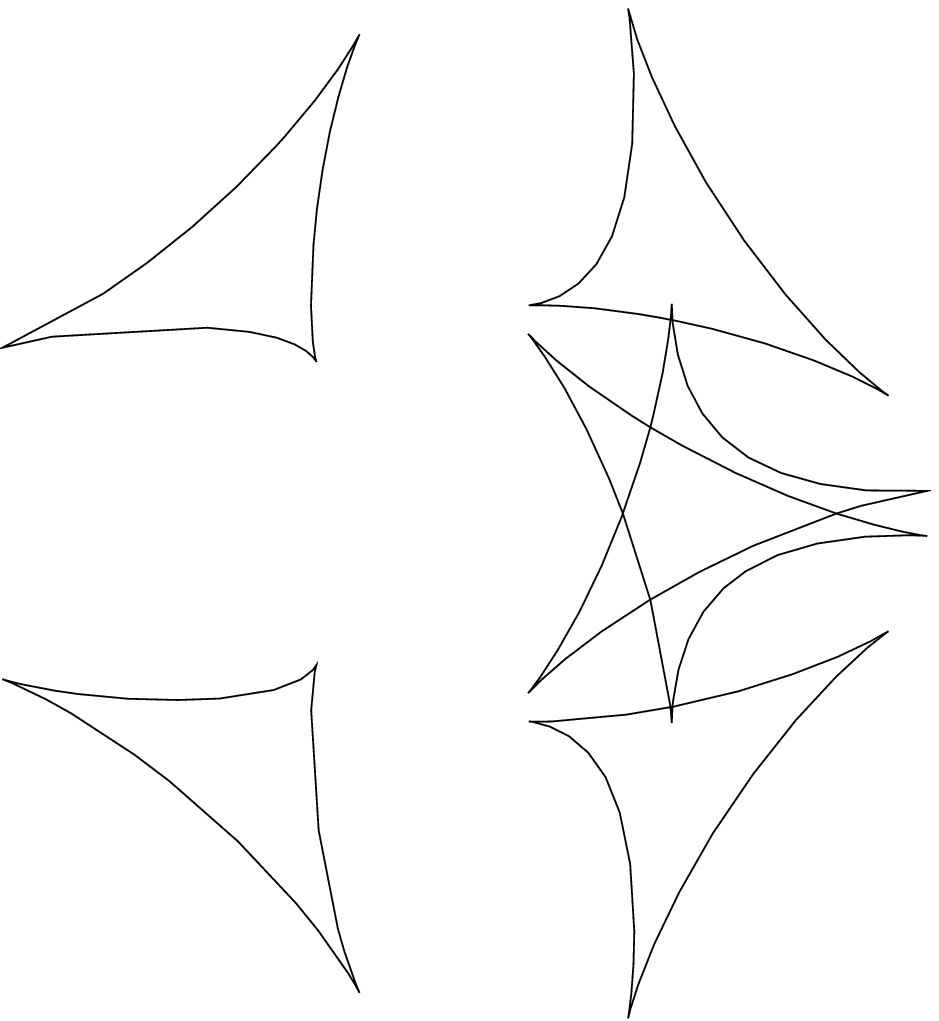}
\end{center}
\caption{\capsize Critical curves of $F_1$ and their images}
\label{fig:criticaldomG}
\end{figure}

Consider the variation on $F_0$
given by $F_1(z) = z^7 + \bar{z}^4 + z$.
Now, the $\bar{z}^4$ term never dominates,
and there is no distinctive intermediate behavior.
The critical set, shown in figure \ref{fig:criticaldomG},
consists of six curves, and from their
images it is clear that each has three outer cusps.
Most points in the image have 7 preimages
and the number of preimages of any point ranges from 7 to 11.
For the program, both functions are in a sense equally easy to handle.

Notice that if one of the critical curves had somehow escaped detection
then the identity from proposition \ref{prop:rotatepolydisk} would have
indicated that something was missing.
More explicitly, we would consider $\Delta_5$, a disk with 5 holes:
$\gamma_0$ would be a large positively oriented circle
and $r(F(\gamma_0)) = 7$ (the degree of $F_1$).
On the other hand, $\gamma_i$ for $i = 1, \ldots, 5$
would be smooth closed curves just outside known critical curves;
from proposition \ref{prop:rotatecritical}, $r(F(\gamma_i)) = 2$.
Here $s = 1$ ($\det DF(p) > 0$ for large $p$) and $n = 5$;
proposition \ref{prop:rotatepolydisk} would give
$7 = 2 + 2 + 2 + 2 + 2 - 4 = 6$;
the fact that this is wrong indicates that at least one critical
curve is missing.

\begin{figure}[ht]
\vglue 11pt
\begin{center}
%\epsfig{height=40mm,file=lip.dom.pln.eps}
%\qquad \qquad
\epsfig{height=40mm,file=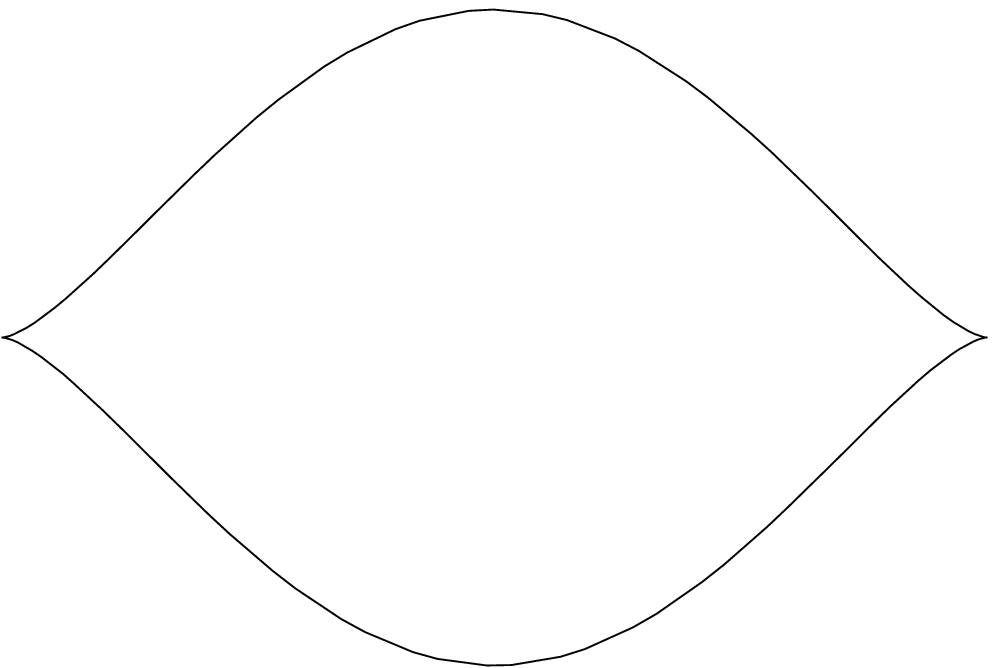}
\end{center}
\caption{\capsize The lip}
\label{fig:lip}
\end{figure}

The {\it lip} is a simpler example:
$F_2(x,y) = (x, y^3/3 + (x^2 - 1)y)$,
with critical set $C$ equal to the unit circle 
and image $F_2(C)$ given in figure \ref{fig:lip}.
% This function has a remarkable property:
% the image of a large circle $\gamma$ centered at the origin
% is a positively oriented simple closed curve.
There exists a diffeomorphism $\tilde F_2$ os the plane
which coincides with $F_2$ outside a certain circle around the origin.
Thus, the propositions in section 9 would not help us detect topological lips.

\begin{figure}[ht]
\vglue 11pt
\begin{center}
\epsfig{height=55mm,file=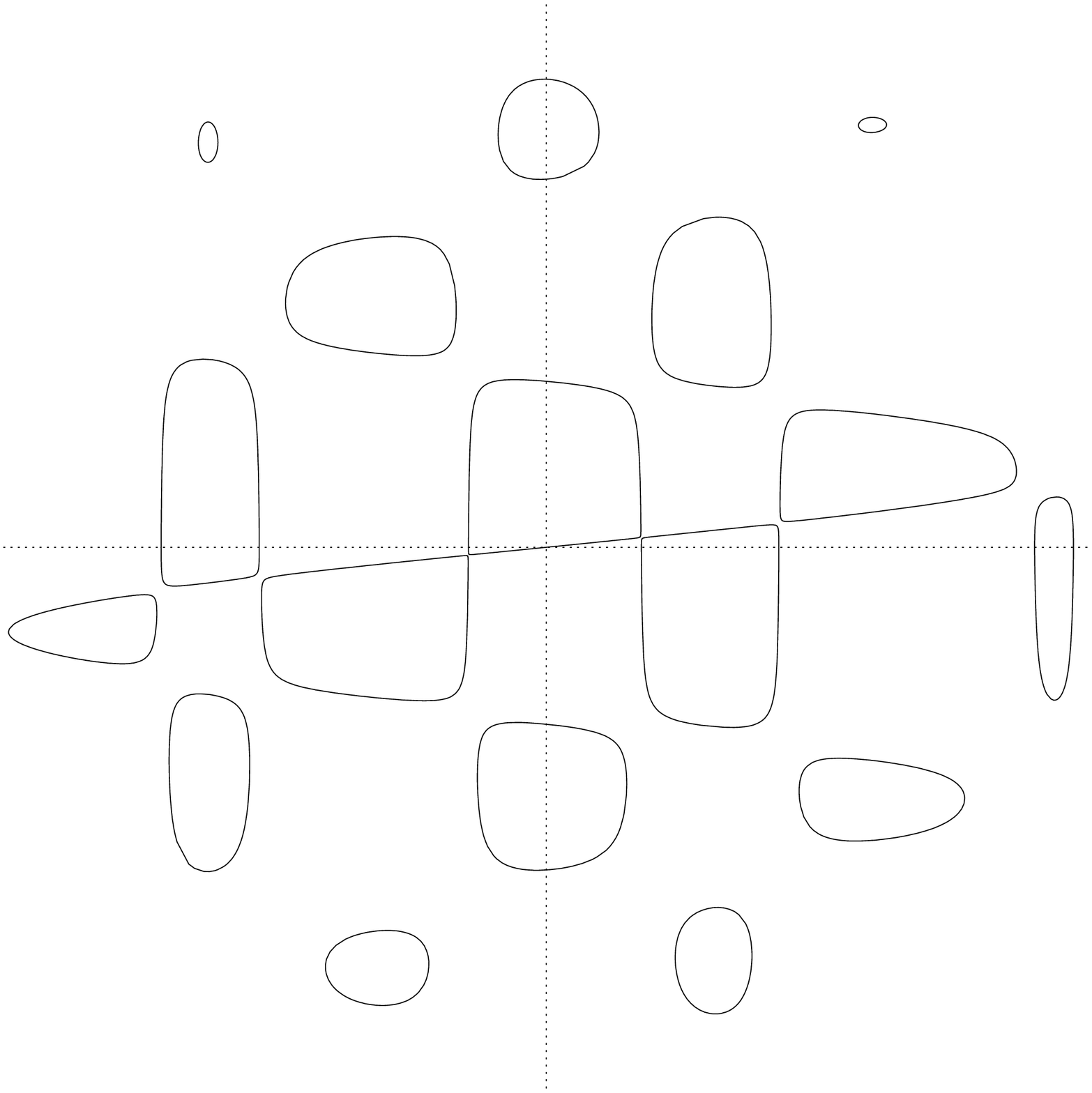}
\qquad \qquad
\epsfig{height=55mm,file=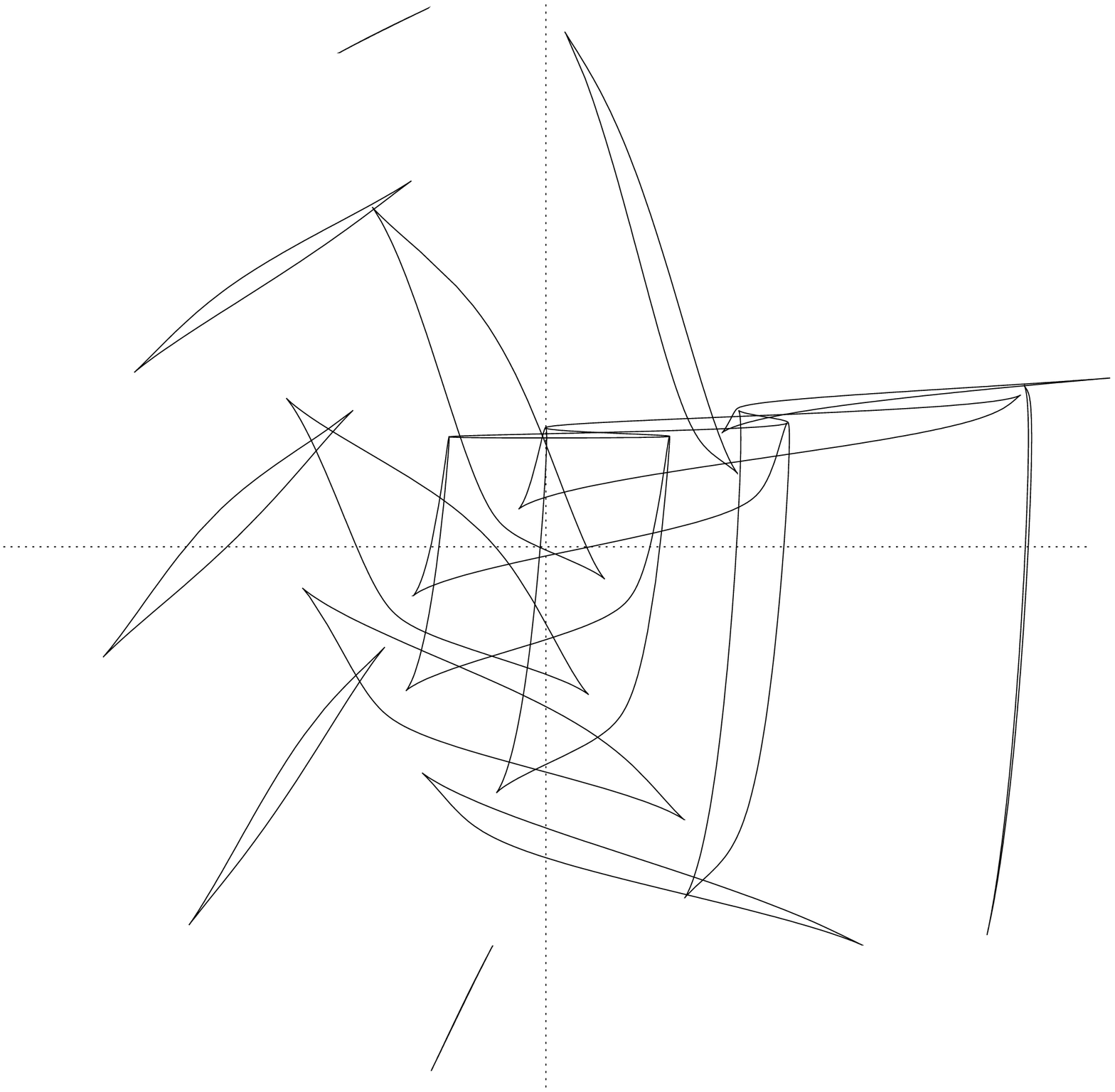}
\end{center}
\caption{\capsize A periodic perturbation of $z \mapsto z^2$}
\label{fig:sincos}
\end{figure}

Our final example is the function
$F_3(x,y) = (x^2 - y^2 + 20 \sin x , 2 x y + 20 \cos y)$.
Both $C$ and $F(C)$ are given in figure \ref{fig:sincos};
the critical set has 17 components and lips abound.
Points in the unbounded tile for $F(C)$ have two preimages;
the origin has 10 preimages.

\bigskip

\vbox{\obeylines

{\nobf Acknowledgements:}
This work was supported by CNPq and Faperj (Brazil).

{\nobf Address:}
Departamento de Matem\'atica, PUC-Rio,
R. Marqu\^es de S. Vicente 225, Rio de Janeiro, RJ 22453-900, Brazil,
{\tt nicolau@mat.puc-rio.br}, {\tt tomei@mat.puc-rio.br}

}

\bigbreak

% \bigskip\bigskip\bigbreak

% {

% \parindent=0pt
% \parskip=0pt
% \obeylines

% Nicolau C. Saldanha and Carlos Tomei, Departamento de Matem\'atica, PUC-Rio
% R. Marqu\^es de S. Vicente 225, Rio de Janeiro, RJ 22453-900, Brazil

% \smallskip

% nicolau@mat.puc-rio.br; http://www.mat.puc-rio.br/$\sim$nicolau/
% tomei@mat.puc-rio.br

% }


\begin{thebibliography}{[10]}

\bibitem{2x2}{\tt http://www.mat.puc-rio.br/$\sim$nicolau/2x2/2x2.html}.
\bibitem{Ahlfors}{ Ahlfors, L. V.,
{Complex Analysis},
McGraw Hill, London, 1981.}
\bibitem{AG}{ Allgower, E.~L. and Georg, K.,
{Numerical continuation methods: an introduction},
Springer-Verlag, New York, 1991.}
\bibitem{BN}{ Bak, J. and Newman, D.~J.,
{Complex analysis},
UTM, Springer-Verlag, New York, 1982.}
\bibitem{CS}{ Chinn, W. G. and Steenrod, N. E.,
{First concepts of topology},
New Mathematical Library, MAA, 1966.}
\bibitem{Duczmal}{ Duczmal, L.,
{Geometria e invers\~ao num\'erica de fun{\c c}\~oes
de uma regi\~ao limitada do plano no plano},
Ph.~D. Thesis, PUC-Rio, Rio de Janeiro, 1997.}
\bibitem{FT}{ Francis,~G.~K. and Troyer,~S.~F.,
{Excellent maps with given folds and cusps},
Houston J. of Math. 3 (1977), 165--192.}
\bibitem{GG}{ Golubitsky,~M. and Guillemin,~V.,
{Stable mappings and their singularities},
Graduate Texts in Mathematics 14,
Springer-Verlag, New York, 1973.}
\bibitem{GS}{ Golubitsky, M. and Schaeffer, D.,
{Singularities and groups in bifurcation theory, vol. 1},
Applied Mathematical Sciences, 51,
Springer-Verlag, New York, 1985.}
\bibitem{H}{ Hopf, H.,
{\"Uber die Drehung der Tangenten und Sehnen ebener Kurven},
Compositio Math. 2 (1935), 50--62.}
\bibitem{L}{ Lang, S.,
{Analysis I},
Addison-Wesley, Reading, MA, 1968.}
\bibitem{Lehto}{ Lehto, O.,
{Univalent Functions and Teichm\"uller Spaces},
GTM 109, Springer-Verlag, New York, 1987.}
\bibitem{MST1}{ Malta, I., Saldanha, N.~C. and  Tomei, C.,
{Critical Sets of Proper Whitney Functions in the Plane},
Matem\'atica Contemporânea, SBM, vol. 13
({10th Brazilian Topology Meeting}), 181--228 (1997).}
\bibitem{MST2}{ Malta, I., Saldanha, N.~C. and  Tomei, C.,
{The numerical inversion of functions from the plane to the plane},
Mathematics of Computation 65, no. 216, 1531--1552 (1996).}
\bibitem{Massey}{ Massey, W.~S.,
{A basic course in algebraic topology},
Graduate Texts in Mathematics 127,
Springer-Verlag, New York, 1991.}
% \bibitem{MT}{Mohar, B. and Thomassen, C.,
% {Graphs on surfaces},
% John Hopkins Univ. Press, Baltimore, MD, 2001.}
\bibitem{Munkres}{ Munkres, J.~R.,
{Topology: a first course},
Prentice-Hall, Inc., Englewood Cliffs, NJ, 1975.}
\bibitem{OR}{ Ortega, J.~M. and Rheinboldt, W.~C.,
{Iterative solution of nonlinear equation in several variables},
Academic Press, New York, 1970.}
\bibitem{Poe}{ Po\'enaru, V.,
{Extending immersions of the circle (d'apr\`es Samuel Blank)},
Expos\'e 342, S\'eminaire Bourbaki 1967-68, Benjamin, NY, 1969.}
\bibitem{Rudin}{ Rudin, W.,
{Real and Complex Analysis},
Third edition, McGraw-Hill, New York, 1987.}
\bibitem{Thomassen}{ Thomassen, C.,
{The Jordan-Sch\"onflies theorem and the classification of surfaces}
Amer. Math. Monthly, 99 , 2, 116--130 (1992). }
\bibitem{Troyer}{ Troyer, S.~F.,
{Extending a boundary immersion to the disk with $n$ holes},
Ph.~D. Thesis, Northeastern U., Boston, MA, 1973.}
\bibitem{W1}{ Whitney, H.,
{On regular closed curves in the plane},
Compositio Math. 4 (1937), 276--284.}
\bibitem{W}{ Whitney, H.,
{On singularities of mappings of Euclidean spaces, I:
mappings of the plane into the plane},
Ann. of Math. 62 (1955), 374--410.}


\end{thebibliography}
\end{document}